\documentclass{article}
\usepackage{amssymb,amsmath,bm,geometry,graphics,graphicx,url,color}
\usepackage{amsfonts}
\usepackage{graphicx}
\usepackage{multirow}

\def\sinc{\mathrm{sinc}}
\def\pa{\partial}
\def\pdt2{\partial_t^2}
\def\pdx2{\partial_x^2}
\newcommand{\norm}[1]{\left\Vert#1\right\Vert}
\newcommand{\normmm}[1]{{\left\vert\kern-0.25ex\left\vert\kern-0.25ex\left\vert #1
    \right\vert\kern-0.25ex\right\vert\kern-0.25ex\right\vert}}
\newcommand{\weig}[1]{\langle#1\rangle}
\newcommand{\abs}[1]{\left\vert#1\right\vert}

\def\TT{{\mathbb{T}}}
\def\RR{{\mathbb{R}}}
\def\ZZ{{\mathbb{Z}}}
\def\ii{\mathrm{i}}

\newtheorem{theo}{Theorem}[section]
\newtheorem{lem}[theo]{Lemma}

\newtheorem{rem}[theo]{Remark}
\newtheorem{defi}[theo]{Definition}

\newtheorem{assum}[theo]{Assumption}
\newtheorem{prop}[theo]{Proposition}

\def\proof{\noindent\underline{Proof}\quad}
\def\QED{\mbox{~$\Box{~}$}}
\def\no{\noindent}

\title{One-stage explicit trigonometric  integrators\\ for effectively solving quasilinear wave equations}

\author{Bin Wang\,\footnote{School of Mathematics and Statistics,
Xi'an Jiaotong University, Xi'an, Shannxi   710049, P.R.China.
E-mail:~{\tt wangbinmaths@qq.com}} \and Changying Liu
\thanks{School of Mathematics and Statistics, Nanjing University of
Information Science and Technology, Nanjing 210044, P.R.China.
E-mail:~{\tt chyliu88@gmail.com}}  \and Yonglei Fang
\thanks{School of Mathematics and Statistics, Zaozhuang
University, Zaozhuang   277160, P.R.China. E-mail:~{\tt
ylfangmath@163.com}}}

\begin{document}
\maketitle

\begin{abstract}
In this paper,   one-stage explicit trigonometric integrators for
solving quasilinear wave equations are formulated and studied. For
solving wave equations, we first introduce trigonometric integrators
as the semidiscretization in time and then consider a spectral
Galerkin method for the discretization in space. We show that
one-stage explicit trigonometric integrators in time have
second-order convergence and the result is also true for the fully
discrete scheme  without requiring any CFL-type coupling of the
discretization parameters. The results are proved by using  energy
techniques, which are widely applied  in the numerical analysis of
methods for partial differential equations.
\medskip

\no{Keywords:} quasilinear wave equations, trigonometric
integrators, second-order convergence, energy technique

\medskip
\no{MSC:} 65M15, 65P10, 65L70, 65M20.

\end{abstract}

\section{Introduction}\label{intro}
In this paper, we are  devoted to the numerical methods  for
effectively solving quasilinear wave equations of the   form (see
\cite{Gauckler17})
\begin{equation}\label{wave equa}
\begin{array}[c]{ll}
\pdt2 u=\pdx2 u-u+\kappa a(u)\pdx2 u+\kappa g(u,\pa_xu), \quad x\in
\TT = \RR/(2\pi\ZZ),\ t\in[0,T],
\end{array}
\end{equation}
 with the smooth and real-valued functions $g$ and $a$  satisfying
$g(0, 0) = a(0) = 0.$ In this paper,   the strength of the
nonlinearities is emphasized by  the   real-valued parameter
$\kappa$  and we  consider $\kappa$  to be small $0<\kappa\ll 1$
such that the nonlinearities are small. The initial values at time
$t = 0$ are assumed to be
\begin{equation}\label{initial val}
\begin{array}[c]{ll}
 u(\cdot,0)=u_0,\ \ \  \pa_tu (\cdot,0)=\dot{u}_0
\end{array}
\end{equation}
and  the boundary conditions are $2\pi$-periodic in one space
dimension. It is noted that the  solutions of \eqref{wave equa} are
assumed to be real-valued in this paper.

It is well known that quasilinear  wave equations occur in   a
variety of applications such as elastodynamics and general
relativity (see, e.g. \cite{Gerner16,Hormander97,Taylor91}).
 These equations have also been used to
 describe many problems which appear in  elasticity,
fluid mechanics and general relativity (see, e.g. \cite{Hughes-77}).
Compared with many publications about the analysis of these
equations
(\cite{Cheng2019,Gerner16,Hormander97,Hughes-77,Taylor91,wang-2019}),
there is much less work devoted to the numerical solutions and
numerical analysis for quasilinear  wave equations.

  These equations with
small $\kappa$ have been extensively studied by
\cite{Dull-15,Chong-13,Dull-17,Groves05}.
 However, in the numerical discretization of
\eqref{wave equa}, the quasilinear term $\kappa a(u)\pdx2 u$ is the
principal difficulty, which needs to be dealt with carefully. In
order to effectively solve \eqref{wave equa},   some implicit and
semi-implicit methods of Runge-Kutta type for semi-discretization in
time   were proposed and researched recently in
\cite{Hochbruck17,Lubich17} . More recently, the authors in
\cite{Gauckler17} showed that a class of explicit exponential
integrators given in \cite{hairer2002,Hochbruck2010} can be used to
numerically solve the quasilinear wave equation \eqref{wave equa}
with two regimes of $\kappa$ by using the energy technique with a
modified discrete energy.

In order to effectively  solve the quasilinear wave equation
\eqref{wave equa},  a class of one-stage explicit trigonometric
integrators will be  rigorously  studied in this paper. We prove
second-order convergence not only for the methods in time  but also
for the fully discrete schemes.
  These trigonometric integrators  were firstly developed   in \cite{wu2010-1} for solving highly oscillatory ODEs and we  refer the reader to
\cite{wang-2016,wang-2018,wu2013-ANM,wang2017-Cal,wu2013-book} for
further researches. Meanwhile, this kind of methods has been applied
to wave equations in the semilinear case (see, e.g.
\cite{Liu_Iserles_Wu(2017-2),JCP-Liu,wang-IMA,JCAM(2016)_Wu_Liu,wubook2015}).
However, these methods have not been researched for quasilinear wave
equations, which motivates this paper.

The  main  contribution of this work is  to show the error bounds of
trigonometric integrators for quasilinear wave equations. In
contrast to the analysis in \cite{Gauckler17}, we do not use a
modified discrete energy  in this paper and just take the simple and
normal energy technique, which is widely used in  the numerical
analysis of partial  differential equations (see, e.g.
\cite{Cano13,Cohen08,Dong14,Gauckler15,Thalhammer15,Thalhammer16,Hochbruck17,Lubich17,Lubich95,wang-2018}).
The  paper is displayed  as follows. In Section
\ref{sec-TRIGONOMETRIC spectral} trigonometric integrators for the
discretization in time  and full-discrete trigonometricintegrators
  are introduced.
The main results of this paper are presented in Section
\ref{sec-main results} and a numerical experiment is carried out  to
show the numerical behaviour and support the theoretical analysis.
In Section \ref{sec-Global error} we   prove  the error bounds for
trigonometric integrators in time. Section \ref{sec-Global error2}
is devoted to the proof of  error bounds for full-discrete
trigonometric integrators.  For one of the trigonometric
integrators, a simple proof for the error bounds  is presented in
Section \ref{sec-simple proof} by establishing a relationship
between this integrator  and a trigonometric integrator researched
in \cite{Gauckler17}.  Finally, in Section \ref{sec-conclu} we
include the conclusions of this paper.

\section{Trigonometric integrators}\label{sec-TRIGONOMETRIC spectral}
In this paper, we will use the following notations and properties,
which have been used in \cite{Gauckler17}.
\begin{itemize}
\item Denote by $H^s = H^s(\TT)$ with $s\geq0$ the usual Sobolev
space and its norm $\norm{\cdot}_s$ is given by
\begin{equation}\label{weighted norm}
\norm{v}_s^2=\sum\limits_{j\in\ZZ}\weig{j}^{2s} |\hat{v}_j|^2\quad
\textmd{for} \quad v(x)=\sum\limits_{j\in\ZZ}
 \hat{v}_j\mathrm{e}^{\mathrm{i}jx},
\end{equation}
where  the weights $\weig{j}$ for $j\in\ZZ$ are defined
$\weig{j}=\sqrt{j^2+1}.$

\item  The corresponding scalar product is
defined by $\langle\cdot\rangle_s$:
\begin{equation*}
\langle v,w\rangle_s=\sum\limits_{j\in\ZZ}\weig{j}^{2s}
 \bar{\hat{v}}_j \hat{w}_j\quad \textmd{for} \quad v(x)=\sum\limits_{j\in\ZZ}
 \hat{v}_j\mathrm{e}^{\mathrm{i}jx},\ \ \quad w(x)=\sum\limits_{j\in\ZZ}
 \hat{w}_j\mathrm{e}^{\mathrm{i}jx}.
\end{equation*}

\item
The solutions $(u(\cdot, t), \pa_tu(\cdot, t))$ of the quasilinear
wave equation \eqref{wave equa} are studied in the spaces $H^{s+1}
\times H^s$ with the norm
$$\normmm{(u,\dot{u})}_{s}=(\norm{u}^2_{s+1}+\norm{\dot{u}}^2_{s})^{1/2}.$$

\item The following classical estimates in Sobolev spaces will be used in
this paper for $s>\frac{1}{2}$ (see Chapter 13 of \cite{Taylor11}):
\begin{equation}\label{norm pro1}
\norm{u v}_0\leq C_s\norm{u}_0\norm{v}_s,\qquad \norm{u v}_s\leq
C_s\norm{u}_s\norm{v}_s.
\end{equation}

\item Another classical estimates for any smooth function $G$ with  $G(0)
= 0$
 are (see Chapter 13 of \cite{Taylor11})
\begin{equation}\label{norm pro2}
\norm{G(u)}_s\leq \Lambda_s(\norm{u}_s)\norm{u}_s,\quad
\norm{G(u)-G(v)}_s\leq \Lambda_s(\norm{u}_s+\norm{v}_s)\norm{u-v}_s,
\end{equation}
where $\Lambda_s(\cdot)$ is a continuous nondecreasing function.

\item It is noted that the norm and the scalar product have the following
connection
\begin{equation}\label{connection}
\norm{u\pm v}_1^2=\norm{u}_1^2+\norm{v}_1^2\pm 2\langle
u,v\rangle_1.\end{equation}
\end{itemize}

\subsection{Methods for the discretization in time} By
letting \begin{equation}\label{fu} f(u)= a(u)\pdx2 u+ g(u,\pa_xu),
\end{equation}
and the linear operator
\begin{equation*}
\Omega=\sqrt{-\pdx2+1},
\end{equation*}
the quasilinear wave equation \eqref{wave equa} becomes
\begin{equation}\label{wave equa2}
\pdt2 u=-\Omega^2 u +\kappa f(u).
\end{equation}

In what follows,   one-stage explicit trigonometric integrators are
used
for the discretization in time of \eqref{wave equa2}. 

\begin{defi}
\label{trigonometric}  (See \cite{wu2010-1}.) For solving
\eqref{wave equa2}, we consider a one-stage  explicit trigonometric
integrator which   is given by
 \begin{equation}
\left\{\begin{array}
[c]{ll}%
u_{n+c_{1}} &
=\phi_{0}(c_{1}^{2}V)u_{n}+hc_{1}\phi_{1}(c_{1}^{2}V)\dot{u}_{n},\\
u_{n+1} & =\phi_{0}(V)u_{n}+h\phi_{1}(V)\dot{u}_{n}+
 \kappa h^{2}\bar{b}_{1}(V)f(u_{n+c_{1}}),\\
\dot{u}_{n+1} &
=-h\Omega^2\phi_{1}(V)u_{n}+\phi_{0}(V)\dot{u}_{n}+\kappa
h\textstyle b_{1}(V)f(u_{n+c_{1}}),
\end{array}\right.
  \label{methods}%
\end{equation}
where  $c_1\in[0,1]$ denotes a real constant and $h$ is the
stepsize. The coefficients  $b_{1}(V)$ and $\bar{b}_{1}(V)$  are
operator-valued functions of $V\equiv h^{2}\Omega^2$, and further we
define
\begin{equation*}
\phi_{0}(V)=\cos(h\Omega),\qquad  \phi_{1}(V)=\textmd{sinc}(h\Omega):=(h\Omega)^{-1}\sin(h\Omega).%
\label{Phi01}%
\end{equation*}
\end{defi}
From the symmetry conditions of trigonometric integrators given in
\cite{wu2013-book}, it follows that  the   integrator
\eqref{methods} is symmetric if and only if
\begin{equation} \label{sym coeffi}c_1=1/2,\ \
\sinc(h\Omega) b_1(V)=(I+\cos(h\Omega))\bar{b}_1(V),
 \end{equation}
 where $I$ is the identical operator.
Under this condition, the trigonometric integrator  \eqref{methods}
can be rewritten as
 \begin{equation}
\left\{\begin{array}
[c]{ll}%
 u_{n+\frac{1}{2}} &
=\cos( \frac{1}{2}h\Omega)u_{n}+\frac{1}{2}h\sinc(\frac{1}{2}h\Omega)\dot{u}_{n},\\
\dot{u}^{+}_{n}&=\dot{u}_{n}+h \kappa\sinc(h\Omega)^{-1}\bar{b}_{1}(V)f(u_{n+\frac{1}{2} }),\\
\left(
  \begin{array}{c}
    \Omega u_{n+1} \\
    \dot{u}^{-}_{n+1} \\
  \end{array}
\right)
 & = \left(
       \begin{array}{cc}
         \cos(h\Omega) &  \sin(h\Omega) \\
         -\sin(h\Omega) &  \cos(h\Omega) \\
       \end{array}
     \right)\left(
  \begin{array}{c}
    \Omega u_{n} \\
    \dot{u}^{+}_{n} \\
  \end{array}
\right)
 ,\\
\dot{u}_{n+1} & =    \dot{u}^{-}_{n+1}+h
\kappa\sinc(h\Omega)^{-1}\bar{b}_{1}(V)f(u_{n+\frac{1}{2}}).
\end{array}\right.
  \label{methods-2}%
\end{equation}
We denote the numerical flow of this integrator by $\varphi_h$,
i.e., $(u_{n+1},\dot{u}_{n+1})=\varphi_h(u_{n},\dot{u}_{n}).$

In this paper,   the one-stage explicit trigonometric integrator
\eqref{methods} is considered under the following assumption.
\begin{assum}\label{thm ass1}
For the coefficient functions of the one-stage explicit
trigonometric integrator \eqref{methods},  we require the symmetry
condition \eqref{sym coeffi} and assume that there exists a constant
$c$ such that {\begin{equation}
\begin{array}{ll}   \label{ass1}
    |\xi \bar{b}_1(\xi^2)|\leq c  ,&\quad|\xi^{2} \bar{b}_1(\xi^2)|\leq c ,\qquad   \qquad \qquad |\bar{b}_1(\xi^2)-\frac{1}{2}\sinc( \frac{1}{2}\xi)|\leq c
   \xi,
 \\   |\xi b_1(\xi^2)|\leq c,&\quad   |b_1(\xi^2)-\cos( \frac{1}{2}\xi)|\leq c
 \xi^{2},\ \ \ \ \  |\xi \sinc(\xi)^{-1} b_1(\xi^2)|\leq c
\end{array}
\end{equation}}
 for all $\xi\geq 0$
\end{assum}

\subsection{Full-discrete methods}
As the full discretization of \eqref{wave equa2}, we consider the
trigonometric integrators for the discretization in time and a
spectral Galerkin method for the discretization in space (see, e.g.
\cite{Gauckler17}).

Denote the space of trigonometric polynomials of degree $K$  by
$$\mathcal{V}^{K}=\Big\{\sum\limits_{j=-K}^K
 \hat{v}_j\mathrm{e}^{\mathrm{i}jx}:   \hat{v}_j\in \mathbb{C} \Big\}$$
and the $L^2$-orthogonal projection onto this ansatz space by
\begin{equation}\label{tri poly}
\mathcal{P}^{K}(v)= \sum\limits_{j=-K}^K
 \hat{v}_j\mathrm{e}^{\mathrm{i}jx}\quad \textmd{for}\quad  v=
 \sum\limits_{j=-\infty}^{\infty}
 \hat{v}_j\mathrm{e}^{\mathrm{i}jx}\in L^2.
 \end{equation}
 Then   the nonlinearity $f(u)$ in the method in
time  \eqref{methods} is considered to be replaced
 by the following new nonlinearity
\begin{equation}\label{rep fu}
\hat{f}^{K}(u)= \mathcal{P}^{K}(f^{K}(u)),
 \end{equation}
where
$$f^{K}(u)=a^{K}(u)\pdx2 u+  g^{K}(u,\pa_xu)\quad \textmd{with}\quad a^{K}=\mathcal{I}^{K}\circ a,\ g^{K}=\mathcal{I}^{K}\circ g.$$
Here the notation $\mathcal{I}^{K}$ is used to describe the
trigonometric interpolation in the space $\mathcal{V}^{K}$.

We are now in the position to present the fully discrete
trigonometric integrator
 \begin{equation}
\begin{array}
[c]{ll}%
u^K_{n+\frac{1}{2}} &
=\cos(\frac{1}{2}h\Omega)u^K_{n}+\frac{1}{2}h\sinc(\frac{1}{2}h\Omega)\dot{u}^K_{n},\\
u^K_{n+1} & =\cos(h\Omega)u^K_{n}+h\sinc(h\Omega)\dot{u}^K_{n}+
 \kappa h^{2}\bar{b}_{1}(V)\hat{f}^K(u^K_{n+\frac{1}{2}}),\\
\dot{u}^K_{n+1} &
=-\Omega\sin(h\Omega)u^K_{n}+\cos(h\Omega)\dot{u}^K_{n}+\kappa
h\textstyle b_{1}(V)\hat{f}^K(u^K_{n+\frac{1}{2}}),
\end{array}
  \label{methods-full}%
\end{equation}
where $u^K_{n}\in \mathcal{V}^{K}$ and $\dot{u}^K_{n}\in
\mathcal{V}^{K}$ are the numerical solutions of   $u(\cdot,t_n)$ and
$\pa_t u(\cdot,t_n)$ respectively. Moreover, the initial values
$u_0$ and $\dot{u}_0$ of \eqref{initial val} are replaced by
$$u^K_{0}=\mathcal{P}^{K}(u_0),\qquad \dot{u}^K_{0}=\mathcal{P}^{K}(\dot{u}_0).$$
We denote the fully discrete  integrator \eqref{methods-full} as
$(u^K_{n+1},\dot{u}^K_{n+1})=\varphi^K_h(u^K_{n},\dot{u}^K_{n}).$

\begin{rem}
It is noted that the nonlinearity $\hat{f}^{K}$ appearing in the
fully discrete trigonometric integrator \eqref{methods-full}  can be
computed efficiently by fast Fourier techniques (see
\cite{Gauckler17}).


\end{rem}

\section{Main results and numerical test}\label{sec-main results}
In this section,   the error bounds are presented not only for the
methods in time but also for the fully discrete schemes.  The exact
solution $u(x, t)$ to \eqref{wave equa2} is required to  satisfy the
following assumption, which has been considered in
\cite{Gauckler17}.
\begin{assum}\label{thm ass2}
(See \cite{Gauckler17}.) The exact solution $(u(\cdot, t),\pa_t
u(\cdot, t))$ to \eqref{wave equa2}  is assumed to be in
$H^{5+s}\times H^{4+s}$ with
\begin{equation} \label{ass2-1}
\normmm{(u(\cdot, t),\pa_t u(\cdot, t))}_{4+s}\leq M \quad
\textmd{for}\quad 0\leq t\leq T,
 \end{equation}
where $s \geq0$ and $M > 0$. Moreover, we assume that there are
 $0 < \delta < 1$ and $A_0 \geq0$
such that $ 1+\kappa a(u(\cdot, t))\leq \delta>0 $ and $
 \kappa a(u(\cdot, t))\leq A_0$
for $0\leq t\leq T.$
\end{assum}
\begin{rem} The regularity assumption
\eqref{ass2-1} on the exact solution was considered in
\cite{Gauckler17}
 and it is true locally in time for initial values in $H^{5+s}\times
H^{4+s}$ by local well-posedness theory (see
\cite{Hughes-77,Taylor91}).
\end{rem}

\subsection{Main results}
\begin{theo}\label{thm-error bounds}
(\textbf{Convergence for trigonometric   integrators in time.})
Assume that Assumption \ref{thm ass1} holds for the coefficient
functions of trigonometric integrators and  Assumption \ref{thm
ass2} is true for the exact solution $(u(\cdot, t),\pa_t u(\cdot,
t))$ with $s = 0$. Then there is a constant $h_0>0$ such that for
$0<\kappa\ll 1$ and for all { $\kappa \lesssim h<h_0$,}
 the following convergence for the time-discrete trigonometric integrator
  $(u_n, \dot{u}_n)$
\eqref{methods}  in $H^2 \times H^1$ holds
\begin{equation}\label{main rel 1}\normmm{(u_n, \dot{u}_n)-(u(\cdot, t_n),\pa_t u(\cdot, t_n))}_{1} \leq Ch^{2} \qquad
\textmd{for} \qquad 0\leq t_n=nh \leq T,\end{equation} where the
constant  $C$
  {depends on  the} smooth functions $a$ and $g$ in
\eqref{wave equa}, the constant $c$ of  Assumption \ref{thm ass1},
the constant  $M$ from Assumption \ref{thm ass2}, but
 is independent of the time step-size $h$ and the final time $T$.
\end{theo}

\begin{theo}\label{thm-error bounds full}
(\textbf{Convergence for full-discrete trigonometric integrators.})
Under the conditions in Theorem \ref{thm-error bounds} but with a
fixed $s\geq 0$ instead of $s=0$,
 there is   $h_0>0$ such that for  $0<\kappa\ll 1$ and for
 all {
 $\kappa
\lesssim h<h_0$,}
 the convergence for the full-discrete trigonometric integrator
$(u^{K}_n, \dot{u}^{K}_n)$ of \eqref{methods-full}  in $H^2 \times
H^1$ is
\begin{equation}\label{main rel 2}\normmm{(u^{K}_n, \dot{u}^{K}_n)-(u(\cdot, t_n),\pa_t u(\cdot, t_n))}_{1} \leq Ch^{2}+CK^{-s-2} \qquad
\textmd{for} \qquad 0\leq t_n=nh \leq T.\end{equation}
\end{theo}

\begin{theo}\label{thm-error bounds a}
(\textbf{Convergence for a special trigonometric integrator.}) For a
special trigonometric integrator TI2 which is presented in Table
\ref{praTRIGONOMETRIC} and does not satisfy  the last requirement in
Assumption
 \ref{thm ass1}, it has the  global error bound \eqref{main rel
1} for TI2 in time and
 the global error bound \eqref{main rel 2} for the full-discrete TI2.
\end{theo}

\begin{rem}
It is noted that this paper only considers one regime of $\kappa$
which is that $\kappa$ is small. The reason is that the bound
\eqref{remainder bound} can be  true only for this case. For the
regime $1$ of $\kappa$, we can only obtain that $$
\abs{\mathcal{R}(u_{n+1},v_{n+1},u_{n},v_{n})}\leq C_M
h\normmm{(u_{n}-v_{n},\dot{u}_{n}-\dot{v}_{n})}_{1},$$ which is not
sufficient for deriving the second-order convergence of the
trigonometric integrators.  For the convergence of the trigonometric
integrators when applied to quasilinear wave equations with
$\kappa=1$,   the possible way to work is to use a modified energy
instead of the normal  energy techniques and we will study it  in
future.

\end{rem}

\subsection{Numerical test}\label{sec-experiment}
As an example,  we present three practical one-stage explicit
trigonometric integrators and their coefficients are
 listed  in Table \ref{praTRIGONOMETRIC}. \renewcommand\arraystretch{1.6}
\begin{table}$$
\begin{array}{|c|c|c|c|c|c|c|c|}
\hline
\text{Methods} &c_1  &\bar{b}_1(V)   &b_1(V)     \\
\hline
\text{TI1} & \frac{1}{2} &\frac{1}{2}  {\phi^3_{1}(V/4)} & {\phi^2_{1}(V/4)\phi_{0}(V/4) }    \cr
\text{TI2} & \frac{1}{2} &\frac{1}{2} \phi_{1}(V)\phi_{1}(V/4) &\phi_{1}(V)\phi_{0}(V/4)    \cr
\text{TI3} & \frac{1}{2} &\frac{1}{2} \phi_{1}(V)\phi_{1}^2(V/4) &\phi_{1}(V)\phi_{1}(V/4)\phi_{0}(V/4)    \cr
 \hline
\end{array}
$$
\caption{Three one-stage explicit trigonometric integrators.}
\label{praTRIGONOMETRIC}
\end{table} It can be checked easily that
 {these three integrators  except TI2 satisfy all the requirements in Assumption
 \ref{thm ass1}. For the convergence of TI2, we will give another proof in Section \ref{sec-simple proof} which does not rely on Assumption
 \ref{thm ass1}. } For comparison, we choose a trigonometric
 integrator (formula (15) with $c=2$ of \cite{Gauckler17}) and
 denote it as NTI.

\begin{figure}
\includegraphics[width=6cm,height=6cm]{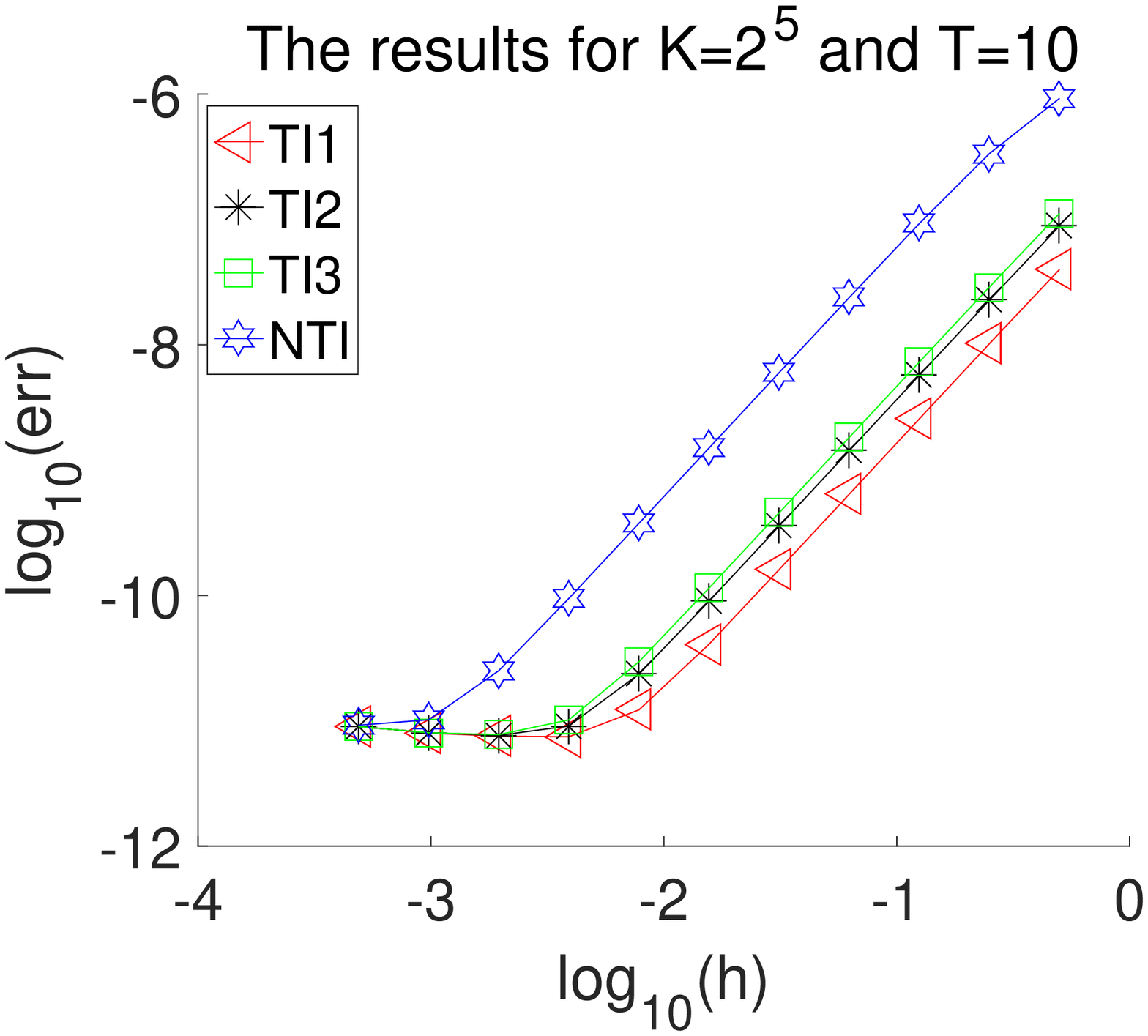}\includegraphics[width=6cm,height=6cm]{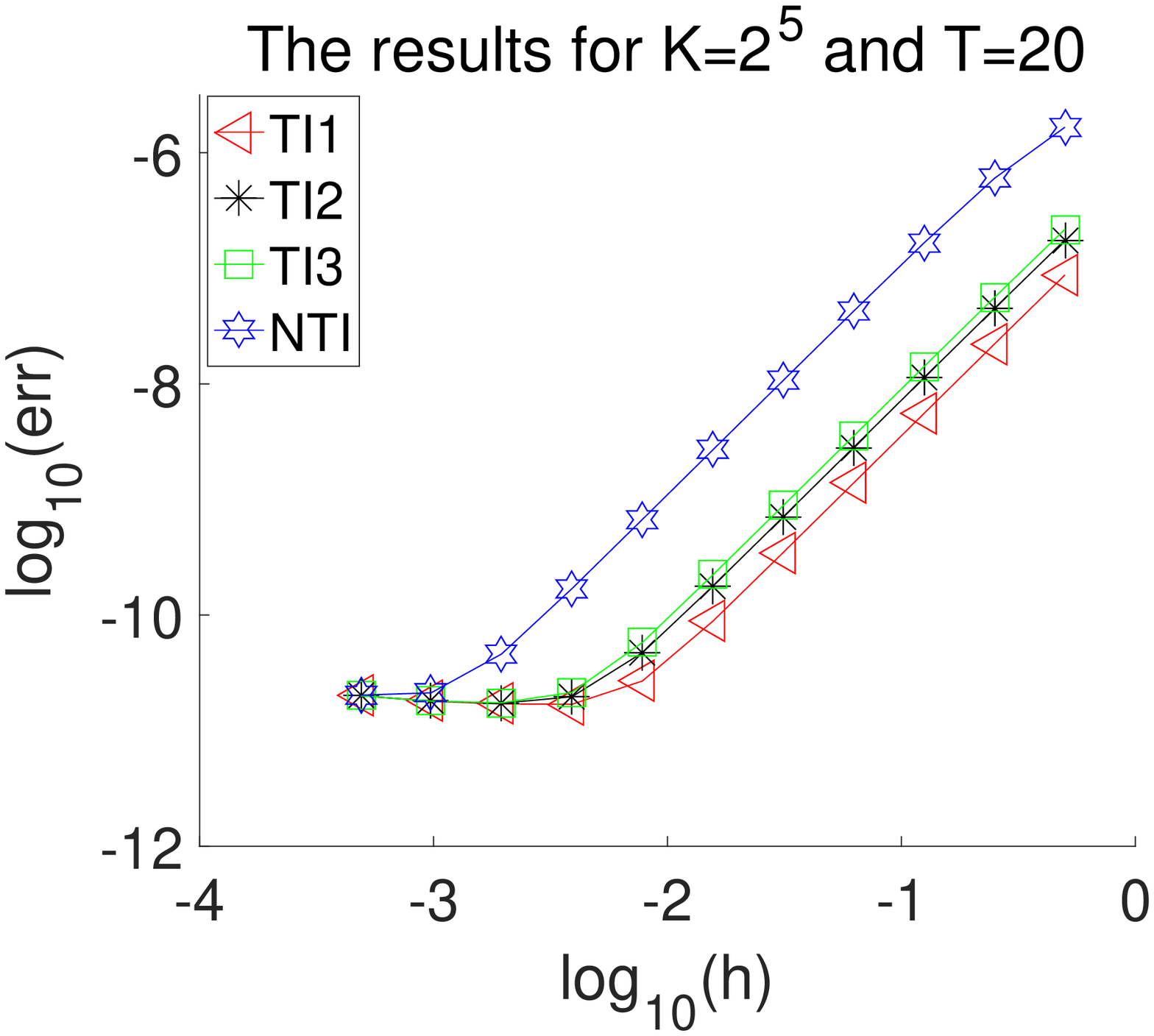}
\caption{Results for $K=2^5$. The logarithm of the  errors against
the logarithm of stepsizes. }
\label{fig1}       
\end{figure}

\begin{figure}
\includegraphics[width=6cm,height=6cm]{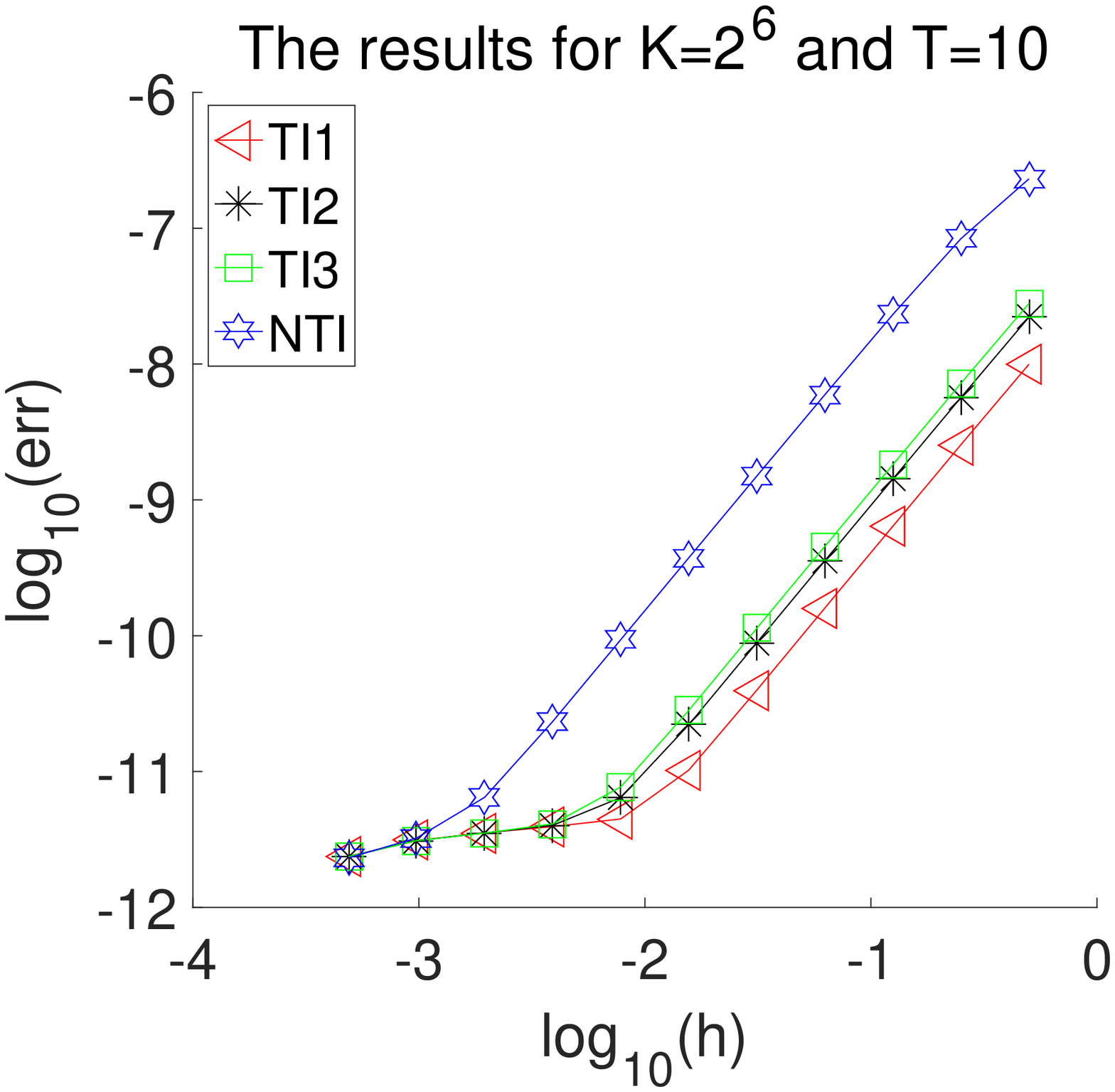}\includegraphics[width=6cm,height=6cm]{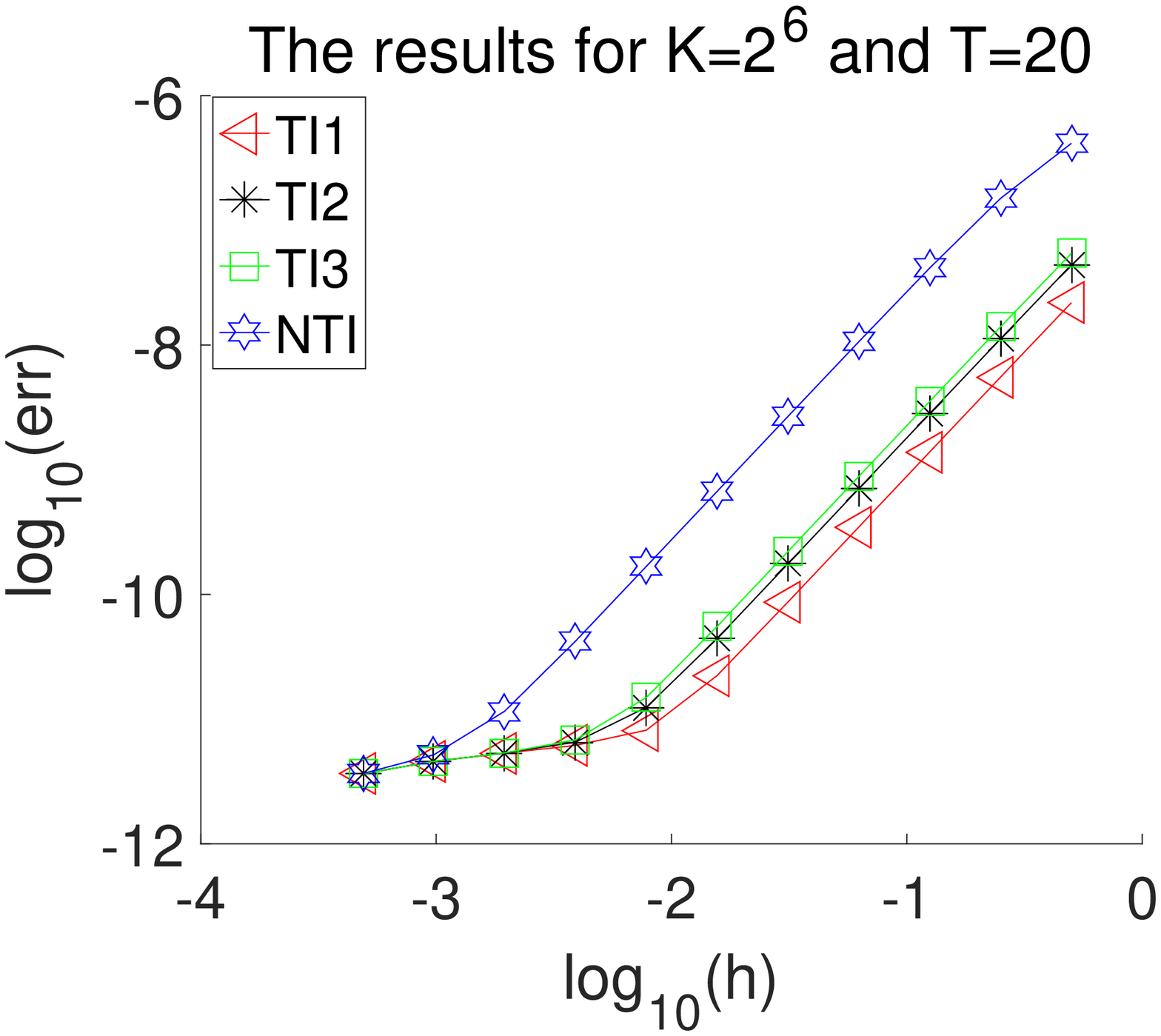}
\caption{Results for $K=2^6$. The logarithm of the  errors against
the logarithm of stepsizes. }
\label{fig2}       
\end{figure}

\begin{figure}
\includegraphics[width=6cm,height=6cm]{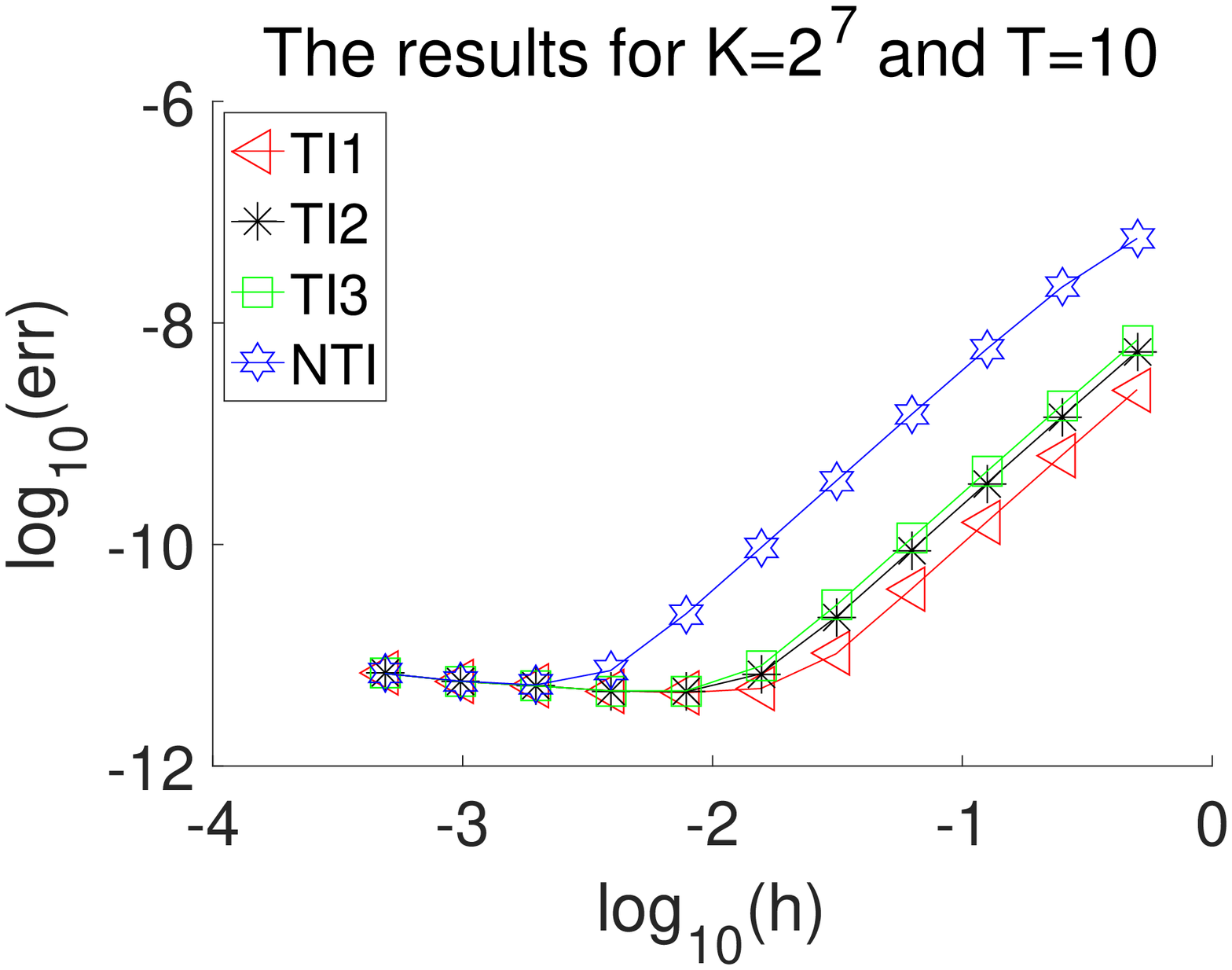}\includegraphics[width=6cm,height=6cm]{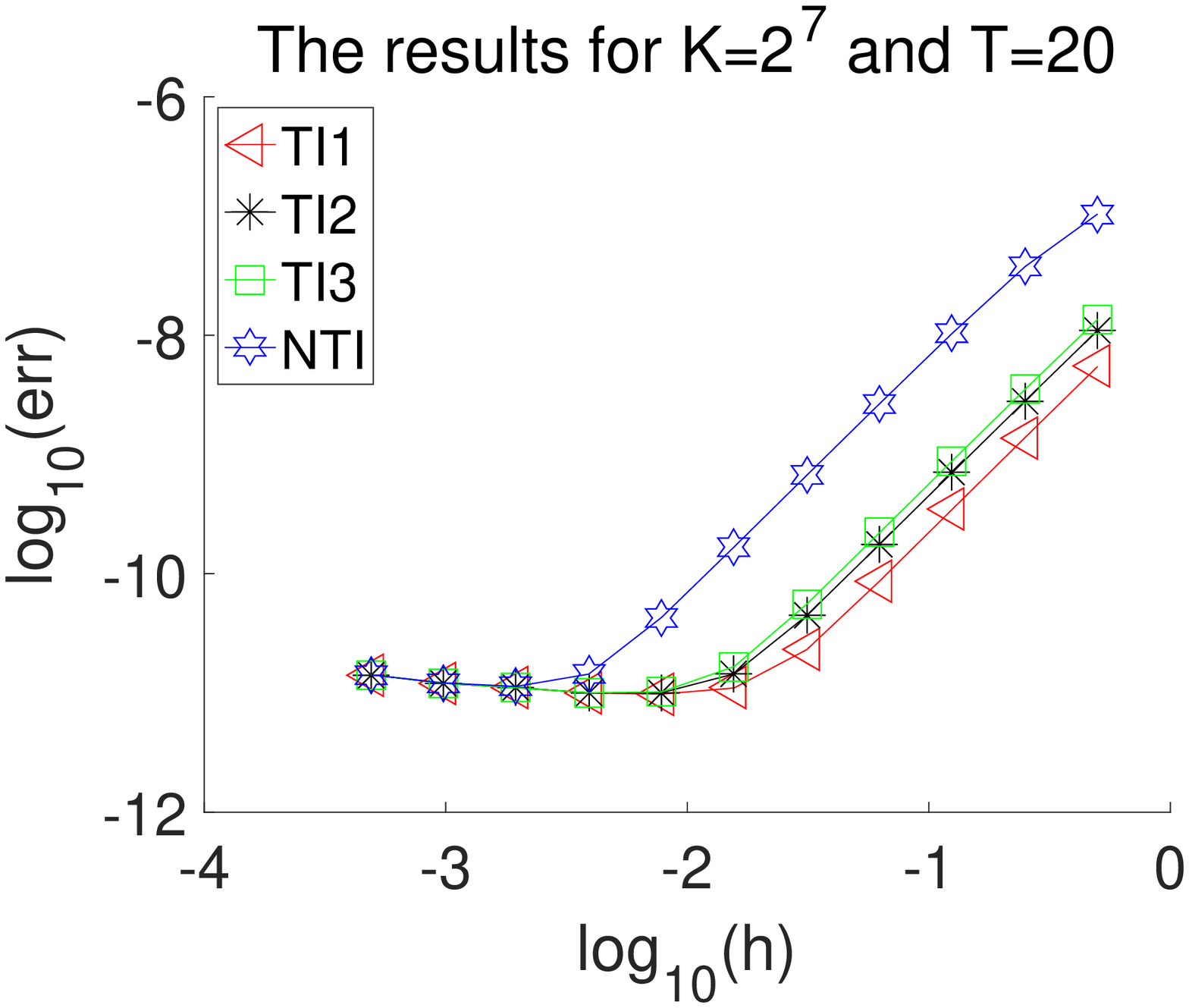}
\caption{Results for $K=2^7$. The logarithm of the  errors against
the logarithm of stepsizes. }
\label{fig3}       
\end{figure}

 We consider the quasilinear wave equation \eqref{wave equa} with
 $a(u)=u$ and $g(u,\pa_xu)=(\pa_xu)^2+\kappa u^3,$
which has been studied in \cite{Chong-13,Gauckler17}. The initial
values are chosen as $$u(x,0)=\sum\limits_{j\in \ZZ}
\frac{1}{\sqrt{1+\abs{j}^{11+\frac{1}{50}}}}e^{\ii jx},\ \ \
\pa_tu(x,0)=\sum\limits_{j\in \ZZ}
\frac{1}{\sqrt{1+\abs{j}^{9+\frac{1}{50}}}}e^{\ii jx}.$$ It is noted
that these initial values are  not in $H^{5+\sigma} \times
H^{4+\sigma}$ for $\sigma\geq 1/100$ but they are in  $H^5 \times
H^4$. Moreover, the regularity assumption \eqref{ass2-1} with  $s =
0$ is true for the initial values. {  We solve this problem  in
$[0,T]$ with $\kappa=1/100$ and the setpsizes $h=\frac{1}{2^j}$ for
$j=1,2,\cdots,11$. The errors in $H^2 \times H^1$ of these three
trigonometric integrators are plotted in Figures
\ref{fig1}-\ref{fig3}.}
 The observed
convergence of these methods are two, which supports the results of
Theorems \ref{thm-error bounds}-\ref{thm-error bounds full}.
Moreover, it follows from the results that the trigonometric
integrators behave better than the
 integrator NTI.

\section{Proof of error bounds for trigonometric integrators in time}\label{sec-Global error}

Theorem \ref{thm-error bounds} will be proved in  this section.
Following \cite{Gauckler17} and in order to present this paper as a
concise proof of concept, we limit ourselves to the exemplary case
$g\equiv0$ in \eqref{wave equa}, i.e.,
\begin{equation} \label{new fu}
f(u)= a(u)\pdx2 u.
 \end{equation}
 Since  the most
critical part of the nonlinearity in  \eqref{wave equa} is the
quasilinear term $a(u)\pdx2 u$, it is straightforward to extend the
proof to nonzero $g$, which will be noted after each step of the
proof.

We remark that in  the proof,  denote by $C$ a generic constant that
may depend on $a$,  the order of the Sobolev space under
consideration and on the constants   in Assumptions \ref{thm ass1}
and \ref{thm ass2}. Denote  by lower indices the additional
dependencies of $C$, e.g., $C_M$ with $M$ from \eqref{ass2-1}.

\subsection{Bounds for a single time step}
By the estimates \eqref{norm pro1}-\eqref{norm pro2} and the
smoothness of $a$, some fundamental properties of the nonlinearity
$f$ in \eqref{new fu} are obtained, which have been given in
\cite{Gauckler17} and will be used in the proof.
\begin{lem}\label{lem-f pro} (See \cite{Gauckler17})
For the nonlinearity $f$ in \eqref{new fu}, it is true that
\begin{equation}\label{fu pro1}
\norm{f(u)}_s\leq
\Lambda_s(\norm{u}_{\sigma})\norm{u}_{\sigma}\norm{u}_{s+2} \quad
\textmd{with}\quad \sigma=\max(s,1),
\end{equation}
and the Lipschitz property
\begin{equation}\label{fu pro2}
\norm{f(u)-f(v)}_s\leq
\Lambda_s(\norm{u}_{s+2}+\norm{v}_{s+2})(\norm{u}_{s+2}+\norm{v}_{s+2})\norm{u-v}_{s+2},
\end{equation}
where   $s \geq0$,  $u, v \in H^{s+2}$, and $\Lambda_s(\cdot)$ is a
continuous non-decreasing function.
\end{lem}

The following lemma shows that the time-discrete trigonometric
integrator $\varphi_h$ given by \eqref{methods-2} maps  $ H^{s+1}
\times H^s$ to itself for $s\geq 1$.

\begin{lem}\label{lem-regularity} (\textbf{Bounds for a single time step.})
Let $s\geq1$ and it is assumed that Assumption \ref{thm ass1} holds.
If  a   time-discrete trigonometric integrator $(u_n,\dot{u}_n)\in
H^{s+1} \times H^s$ satisfies $\normmm{(u_n,\dot{u}_n)}_{s}\leq M,$
then it is true that $$\norm{ u_{n+\frac{1}{2}}}_{s+1}  \leq C_M$$
and $(u_{n+1},\dot{u}_{n+1})\in H^{s+1} \times H^s$ with
$$\normmm{(u_{n+1},\dot{u}_{n+1})}_{s}  \leq C_M.$$
\end{lem}
\proof
 From  the definition of the trigonometric integrator \eqref{methods},
 it follows that
 \begin{equation*}\begin{array}
[c]{ll} \norm{u_{n+\frac{1}{2}}}_{s+1} &\leq
\norm{\cos\big(\frac{1}{2}h\Omega\big)u_{n}}_{s+1}+
\frac{1}{2}\norm{\Omega^{-1}
\sin\big(\frac{1}{2}h\Omega\big)\dot{u}_{n}}_{s+1}\\
&=\norm{\cos\big(\frac{1}{2}h\Omega\big)u_{n}}_{s+1}+
\frac{1}{2}\norm{
\sin\big(\frac{1}{2}h\Omega\big)\dot{u}_{n}}_{s}\leq C_M.
\end{array}
\end{equation*}
Thus
\begin{equation*}\begin{array}[c]{ll}h^2\norm{\bar{b}_1( V)f(u_{n+\frac{1}{2}})}_{s+1}& \leq
  \norm{ \Omega^{-2}f(u_{n+\frac{1}{2}})}_{s+1}= \norm
 {f(u_{n+\frac{1}{2}})}_{s-1}\\
&\leq\Lambda_{s-1}\big(\norm{u_{n+\frac{1}{2}}}_{s+1}\big)
\norm{u_{n+\frac{1}{2}}}^2_{s+1}
\end{array}\end{equation*}
is seen from the second formula in \eqref{ass1} and  \eqref{fu
pro1}. In a similar way, by the fourth formula in \eqref{ass1} it
arrives that
\begin{equation*}\begin{array}[c]{ll}h\norm{b_1( h\Omega)f(u_{n+\frac{1}{2}})}_{s} &\leq
 \norm{ \Omega^{-1}f(u_{n+\frac{1}{2}})}_{s}=\norm
 {f(u_{n+\frac{1}{2}})}_{s-1}\\
&\leq\Lambda_{s-1}\big(\norm{u_{n+\frac{1}{2}}}_{s+1}\big)
\norm{u_{n+\frac{1}{2}}}^2_{s+1}.
\end{array}\end{equation*}
Therefore, considering the scheme of trigonometric integrator
\eqref{methods} again leads to
\begin{equation*}\begin{array}[c]{ll}\norm{ u_{n+1}}_{s+1} &\leq \norm{\cos( h\Omega)u_n}_{s+1}+
  \norm{\Omega^{-1}\sin( h\Omega)\dot{u}_n}_{s+1}+h^2\norm{\bar{b}_1(
h\Omega)f(u_{n+ \frac{1}{2}})}_{s+1}\\
&\leq \norm{ u_n}_{s+1}+
  \norm{ \dot{u}_n}_{s}+\Lambda_{s-1}\big(\norm{u_{n+\frac{1}{2}}}_{s+1}\big) \norm{u_{n+\frac{1}{2}}}^2_{s+1}\end{array}\end{equation*} and
\begin{equation*}\begin{array}[c]{ll}
\norm{ \dot{u}_{n+1}}_{s} &\leq \norm{\Omega\sin( h\Omega)u_n}_{s}+
  \norm{ \cos( h\Omega)\dot{u}_n}_{s}+h\norm{b_1(
h\Omega)f(u_{n+ \frac{1}{2}})}_{s}\\
&\leq \norm{u_n}_{s+1}+
  \norm{  \dot{u}_n}_{s}+ \Lambda_{s-1}\big(\norm{u_{n+\frac{1}{2}}}_{s+1}\big) \norm{u_{n+\frac{1}{2}}}^2_{s+1}.\end{array}\end{equation*}
 \QED

\begin{rem} It is noted that from the proof, it follows that this lemma is still true for  a nonzero $g$ in   \eqref{wave equa}.
\end{rem}

\subsection{Stability}
In this subsection, we will show the stability of trigonometric
integrators. Before presenting the result, the following two lemmas
are needed.

\begin{lem}\label{lem-remainder}
Assume  that Assumption \ref{thm ass1} holds
   with constant $c$. For two time-discrete  trigonometric  numerical solutions $(u_n,\dot{u}_n)\in H^{s+1}
\times H^s$ and $(v_n,\dot{v}_n)\in H^{s+1} \times H^s$ with
$s\geq0$, one has that
\begin{equation*}\begin{array}[c]{ll}
\normmm{(u_{n+1}-v_{n+1},\dot{u}_{n+1}-\dot{v}_{n+1})}^2_{1}
=\normmm{(u_{n}-v_{n},\dot{u}_{n}-\dot{v}_{n})}^2_{1}+\kappa
\mathcal{R}(u_{n+1},v_{n+1},u_{n},v_{n}),\end{array}\end{equation*}
where the remainder is given by
\begin{equation}\begin{array}[c]{ll}\label{remainder}
\mathcal{R}(u_{n+1},v_{n+1},u_{n},v_{n}) \\
=  \langle
2\sinc(h\Omega)^{-1}{b_{1}(V)}(u_{n+1}-u_{n}-v_{n+1}+v_{n}),f(u_{n+\frac{1}{2}})-f(v_{n+\frac{1}{2}})\rangle_1.\end{array}\end{equation}
\end{lem}
\proof In this proof, we will use the following results
\begin{equation}\begin{array}[c]{ll}\label{dun1}
h\sinc(h\Omega)\dot{u}_{n+1}=\cos(h\Omega)u_{n+1}+\kappa h^2
\bar{b}_{1}(V)f(u_{n+\frac{1}{2}
})-u_{n},\\
h\sinc(h\Omega)\dot{u}_{n}=-\cos(h\Omega)u_{n}-\kappa h^2
\bar{b}_{1}(V)f(u_{n+\frac{1}{2}
})+u_{n+1},\end{array}\end{equation} which are obtained by
considering the trigonometric scheme \eqref{methods-2} and its
symmetry. The same relations hold for $v$.

According to   the third step of the integrator \eqref{methods-2},
it is obtained that
\begin{equation}\begin{array}[c]{ll}\label{sta-for1}
\norm{\Omega(u_{n+1}-v_{n+1})}^2_{1}+\norm{\dot{u}^{-}_{n+1}-\dot{v}^{-}_{n+1}}^2_{1}
=\norm{\Omega(u_{n}-v_{n})}^2_{1}+\norm{\dot{u}^{+}_{n}-\dot{v}^{+}_{n}}^2_{1}.\end{array}\end{equation}
By \eqref{connection} and the fourth   step of  \eqref{methods-2},
we have
\begin{equation*}\begin{array}[c]{ll}
&\norm{\dot{u}^{-}_{n+1}-\dot{v}^{-}_{n+1}}^2_{1}=\norm{\dot{u}_{n+1}-\dot{v}_{n+1}
-h
\kappa\sinc(h\Omega)^{-1}\bar{b}_{1}(V)(f(u_{n+\frac{1}{2}})-f(v_{n+\frac{1}{2}}))}^2_{1}\\
=&\norm{\dot{u}_{n+1}-\dot{v}_{n+1}}^2_{1}+h^2
\kappa^2\norm{\sinc(h\Omega)^{-1}\bar{b}_{1}(V)(f(u_{n+\frac{1}{2}})-f(v_{n+\frac{1}{2}}))}^2_{1}\\
&-2h\kappa\langle
\dot{u}_{n+1}-\dot{v}_{n+1},\sinc(h\Omega)^{-1}\bar{b}_{1}(V)(f(u_{n+\frac{1}{2}})-f(v_{n+\frac{1}{2}}))\rangle_1.\end{array}\end{equation*}
Replacing the  difference $\dot{u}_{n+1}-\dot{v}_{n+1}$ with the
help of the first relation of \eqref{dun1} yields
\begin{equation*}\begin{array}[c]{ll}
&h\langle
\dot{u}_{n+1}-\dot{v}_{n+1},\sinc(h\Omega)^{-1}\bar{b}_{1}(V)f(u_{n+\frac{1}{2}})-f(v_{n+\frac{1}{2}})\rangle_1\\
=&\langle \sinc(h\Omega)^{-1}\cos(h\Omega)(u_{n+1}-v_{n+1})+\kappa
h^2\sinc(h\Omega)^{-1} \bar{b}_{1}(V)(f(u_{n+\frac{1}{2}
})-f(v_{n+\frac{1}{2}
}))\\
&-\sinc(h\Omega)^{-1}(u_{n}-v_{n}),\sinc(h\Omega)^{-1}\bar{b}_{1}(V)(f(u_{n+\frac{1}{2}})-f(v_{n+\frac{1}{2}}))\rangle_1\\
=&\langle
\sinc(h\Omega)^{-2}\cos(h\Omega)\bar{b}_{1}(V)(u_{n+1}-v_{n+1}),
f(u_{n+\frac{1}{2}})-f(v_{n+\frac{1}{2}})\rangle_1\\
&+\kappa h^2\norm{\sinc(h\Omega)^{-1}
\bar{b}_{1}(V)(f(u_{n+\frac{1}{2} })-f(v_{n+\frac{1}{2}
}))}_1^2\\
&-\langle \sinc(h\Omega)^{-2}\bar{b}_{1}(V)(u_{n}-v_{n}),
f(u_{n+\frac{1}{2}})-f(v_{n+\frac{1}{2}})\rangle_1.\end{array}\end{equation*}
Here, we use  the property $\langle v,\Psi(h\Omega)
w\rangle_1=\langle \Psi(h\Omega) v, w\rangle_1,$ which is obtained
by Parseval's theorem.

Similarly, taking the second   step of  \eqref{methods-2} and the
second of \eqref{dun1}  into account, one gets
\begin{equation*}\begin{array}[c]{ll}
&\norm{\dot{u}^{+}_{n}-\dot{v}^{+}_{n}}^2_{1}=\norm{\dot{u}_{n}-\dot{v}_{n}}^2_{1}+h^2
\kappa^2\norm{\sinc(h\Omega)^{-1}\bar{b}_{1}(V)(f(u_{n+\frac{1}{2}})-f(v_{n+\frac{1}{2}}))}^2_{1}\\
&+2h\kappa\langle
\dot{u}_{n}-\dot{v}_{n},\sinc(h\Omega)^{-1}\bar{b}_{1}(V)(f(u_{n+\frac{1}{2}})-f(v_{n+\frac{1}{2}}))\rangle_1\end{array}\end{equation*}
and
\begin{equation*}\begin{array}[c]{ll}
&h\langle
\dot{u}_{n}-\dot{v}_{n},\sinc(h\Omega)^{-1}\bar{b}_{1}(V)(f(u_{n+\frac{1}{2}})-f(v_{n+\frac{1}{2}}))\rangle_1\\
=&\langle
-\sinc(h\Omega)^{-2}\cos(h\Omega)\bar{b}_{1}(V)(u_{n}-v_{n}),
f(u_{n+\frac{1}{2}})-f(v_{n+\frac{1}{2}})\rangle_1\\
&-\kappa h^2\norm{\sinc(h\Omega)^{-1}
\bar{b}_{1}(V)(f(u_{n+\frac{1}{2} })-f(v_{n+\frac{1}{2}
}))}_1^2\\
&+\langle \sinc(h\Omega)^{-2}\bar{b}_{1}(V)(u_{n+1}-v_{n+1}),
f(u_{n+\frac{1}{2}})-f(v_{n+\frac{1}{2}})\rangle_1.\end{array}\end{equation*}
In the light of the above analysis, the formula \eqref{sta-for1} can
be expressed as
\begin{equation}\begin{array}[c]{ll}\label{sta-for2}
&\normmm{(u_{n+1}-v_{n+1},\dot{u}_{n+1}-\dot{v}_{n+1})}^2_{1}
-\normmm{(u_{n}-v_{n},\dot{u}_{n}-\dot{v}_{n})}^2_{1}\\
=&2\kappa \langle
\sinc(h\Omega)^{-2}(I+\cos(h\Omega))\bar{b}_{1}(V)(u_{n+1}-u_{n}-v_{n+1}+v_{n}),
f(u_{n+\frac{1}{2}})-f(v_{n+\frac{1}{2}})\rangle_1.\end{array}\end{equation}
From the symmetry condition \eqref{sym coeffi}, it follows that
$$\sinc(h\Omega)^{-2}(I+\cos(h\Omega))\bar{b}_{1}(V)=\sinc(h\Omega)^{-1}b_{1}(V).$$
Therefore, \eqref{sta-for2} yields the statement of this lemma with
the remainder \eqref{remainder}.
 \QED

The bound of the remainder $\mathcal{R}$ \eqref{remainder} is
  estimated by the following lemma.

\begin{lem}\label{lem-remainder bound}
(\textbf{Bound of the remainder.}) Under the conditions given in
Assumption \ref{thm ass1}, if time-discrete  trigonometric numerical
solutions $(u_n,\dot{u}_n)$ and $(v_n,\dot{v}_n)$ belonging to $
H^{3} \times H^2$ satisfy
$$\normmm{(u_{n},\dot{u}_{n})}_{2}\leq M\quad \textmd{and} \quad \normmm{(v_{n},\dot{v}_{n})}_{2}\leq M,$$
we then obtain    the bound for the remainder $\mathcal{R}$ as
\begin{equation}\begin{array}[c]{ll}\label{remainder bound}
\abs{{\kappa}\mathcal{R}(u_{n+1},v_{n+1},u_{n},v_{n})}\leq C_M
h\normmm{(u_{n}-v_{n},\dot{u}_{n}-\dot{v}_{n})}^2_{1}.\end{array}\end{equation}
\end{lem}
 \proof
It is obtained from   the scheme of trigonometric integrators
\eqref{methods} that
 \begin{equation}\label{differ-un}\begin{array}[c]{ll}&\norm{(u_{n+1}-u_{n})-(v_{n+1}-v_{n})}_{1} \\
 \leq &\norm{(\cos( h\Omega)-I)(u_n-v_n)}_{1}+
  \norm{h\sinc( h\Omega)(\dot{u}_n-\dot{v}_n)}_{1}\\
  &+h^2\norm{\bar{b}_1(
h\Omega)(f(u_{n+ \frac{1}{2}})-f(v_{n+ \frac{1}{2}}))}_{1}\\
\leq &2 \norm{\sin(h\Omega/2)^2(u_n-v_n)}_{1}+h
  \norm{ (\dot{u}_n-\dot{v}_n)}_{1}+h\norm{\Omega^{-1}(f(u_{n+ \frac{1}{2}})-f(v_{n+ \frac{1}{2}}))}_{1}\\
  \leq &  \norm{ h\Omega  (u_n-v_n)}_{1}+h
  \norm{(\dot{u}_n-\dot{v}_n)}_{1}+h\norm{f(u_{n+ \frac{1}{2}})-f(v_{n+ \frac{1}{2}})}_{0}\\
 \leq &h\norm{u_n-v_n}_{2}+h
  \norm{\dot{u}_n-\dot{v}_n}_{1}+  h\Lambda_0\Big(\norm{u_{n+ \frac{1}{2}}}_{2}+\norm{v_{n+
  \frac{1}{2}}}_{2}\Big)\\
  &
 \Big( \norm{u_{n+ \frac{1}{2}}}_{2}+\norm{v_{n+ \frac{1}{2}}}_{2}\Big)
 \norm{u_{n+ \frac{1}{2}}-v_{n+ \frac{1}{2}}}_{2},
\end{array}\end{equation}
where the first formula in \eqref{ass1} and \eqref{fu pro2} were
used here.  By Lemma \ref{lem-regularity}, we know that
$$\norm{u_{n+ \frac{1}{2}}}_{2}\leq C_M,\ \ \norm{v_{n+ \frac{1}{2}}}_{2}\leq C_M.$$
Using the scheme of trigonometric integrators \eqref{methods} again
leads
 \begin{equation}\label{rv-}\begin{array}[c]{ll}\norm{u_{n+ \frac{1}{2}}-v_{n+ \frac{1}{2}}}_{2}&\leq \norm{u_{n}-v_{n}}_{2}+ \norm{\sin(
 h\Omega)\Omega^{-1}(\dot{u}_n-\dot{v}_n)}_{2}\\
 &\leq \norm{u_{n}-v_{n}}_{2}+ \norm{\dot{u}_n-\dot{v}_n}_{1}.
\end{array}\end{equation}
By the above results, \eqref{differ-un} becomes
 \begin{equation}\label{differ-uvn}\begin{array}[c]{ll}
 & \norm{(u_{n+1}-u_{n})-(v_{n+1}-v_{n})}_{1}\leq C_M
h\normmm{(u_{n}-v_{n},\dot{u}_{n}-\dot{v}_{n})}_1.
 \end{array}\end{equation}

{Using the facts that $\abs{\langle \cdot, \cdot\rangle_1}\leq C
\norm{\cdot}_1\norm{\cdot}_1$,
 the remainder \eqref{remainder} have the following bound
\begin{equation*}\begin{aligned}
&\abs{\kappa\mathcal{R}(u_{n+1},v_{n+1},u_{n},v_{n})} \\
=&\langle
(u_{n+1}-u_{n}-v_{n+1}+v_{n}),2\sinc(h\Omega)^{-1}b_{1}(V)\kappa(f(u_{n+\frac{1}{2}})-f(v_{n+\frac{1}{2}}))\rangle_1 \\
\leq& C \norm{\Upsilon}_1\norm{2\sinc(h\Omega)^{-1}b_{1}(V)\kappa (
f(u_{n+\frac{1}{2}})-f(v_{n+\frac{1}{2}}))}_1\\
\leq& C \norm{\Upsilon}_1\norm{h(h\Omega)^{-1}(
f(u_{n+\frac{1}{2}})-f(v_{n+\frac{1}{2}}))}_1\\
\leq& C \norm{\Upsilon}_1\norm{
f(u_{n+\frac{1}{2}})-f(v_{n+\frac{1}{2}})}_0,
\end{aligned}\end{equation*}
 where
 $\Upsilon=u_{n+1}-u_{n}-v_{n+1}+v_{n}.$
Based on this fact, the results \eqref{fu pro2},  \eqref{rv-} and
$\norm{\Upsilon}_1\leq C_M
h\normmm{(u_{n}-v_{n},\dot{u}_{n}-\dot{v}_{n})}_1,$  the bound
\eqref{remainder bound} is obtained.}
 \QED

By the above two lemmas, we obtain the following estimate about the
stability of trigonometric integrators.
\begin{prop}\label{pro-Stability}
(\textbf{Stability.}) Under the conditions given in Lemma
\ref{lem-remainder bound}, if the solution  $(u,\pa_t u)$  to
\eqref{wave equa2} in $H^{3} \times H^2$ satisfies
$$\normmm{(u(\cdot,t_n),\pa_t
u(\cdot,t_n))}_2\leq M,$$
then it holds that
\begin{equation}\begin{array}[c]{ll}\label{Stability}
\normmm{(u_{n+1}, \dot{u}_{n+1})-\varphi_h(u(\cdot,t_n),\pa_t
u(\cdot,t_n))}_1^2\leq {(1+C_M h )} \normmm{(e_n,\dot{e}_n)}_1^2,
\end{array}\end{equation}
where the global error $(e_n,\dot{e}_n)$ is defined by
$$(e_n,\dot{e}_n)=(u_{n}, \dot{u}_{n})-(u(\cdot,t_n),\pa_t
u(\cdot,t_n)).$$

\end{prop}
 \proof This result is obtained immediately by letting
$$(v_{n}, \dot{v}_{n}) = (u(\cdot,t_n),\pa_t
u(\cdot,t_n)),\ \ (v_{n+1},
\dot{v}_{n+1})=\varphi_h(u(\cdot,t_n),\pa_t u(\cdot,t_n))$$
 in Lemmas \ref{lem-remainder} and \ref{lem-remainder bound}. \QED

\begin{rem}It is remarked that  the statement of Lemma \ref{lem-remainder} remains valid for a nonzero $g$ in \eqref{wave equa}
with a new remainder $\mathcal{R}^*$ which has additional terms with
$g(u)$ instead of $a(u)\pdx2 u$. The remainder bound given in Lemma
\ref{lem-remainder bound} can be extended to this case since $g(u)$
is more regular than $a(u)\pdx2 u$. Thus the result about the
stability proposed in Proposition \ref{pro-Stability} is still true
for the case that $g$ is nonzero.
\end{rem}

\subsection{Local error bound}
Local error bound of time-discrete  trigonometric integrators is
given by the following proposition.

\begin{prop}\label{prop-local error}
(\textbf{Local error in $H^{2}\times H^{1}$.}) Under the conditions
of Lemma \ref{lem-regularity}, if the solution $(u(\cdot, t),\pa_t
u(\cdot, t))$ to \eqref{wave equa2}
 is  in $H^{5} \times H^4$ with
$$\normmm{(u(\cdot, t),\pa_t u(\cdot, t))}_{4}\leq M,$$
then  one gets
$$\normmm{(d_{n+1},\dot{d}_{n+1})}_{1}\leq {C_Mh^{3}},$$
where the local error $(d_{n+1},\dot{d}_{n+1})$ is defined by
$$(d_{n+1},\dot{d}_{n+1})=\varphi_h(u(\cdot,t_n),\pa_t
u(\cdot,t_n))-(u(\cdot,t_{n+1}),\pa_t u(\cdot,t_{n+1})).$$
\end{prop}
 \proof Without loss of generality, the proof is given in the case $n = 0$, that is, we
 prove the result for
$$(d_{1},\dot{d}_{1})=(u_1,\dot{u}_1)-(u(\cdot,h),\pa_t
u(\cdot,h)).$$ By the variation-of-constants formula, the exact
solution of  \eqref{wave equa2} at $t=h$ can be expressed by
\begin{equation}\label{vari cons}
\begin{array}[c]{ll}
\left(
  \begin{array}{c}
    u(\cdot,h) \\
     \pa_t
u(\cdot,h) \\
  \end{array}
\right)=R(h)\left(
  \begin{array}{c}
   u_0 \\
     \dot{u}_0 \\
  \end{array}
\right)+\kappa\int_{0}^{h}R(h-t)\left(
                                 \begin{array}{c}
                                   0 \\
                                   f(u(\cdot,t)) \\
                                 \end{array}
                               \right)d t,
\end{array}
\end{equation}
where
 \begin{equation}\label{matr R}%
R(t)=\left(
       \begin{array}{cc}
         \cos(t\Omega) &  t \textmd{sinc}(t\Omega) \\
         -\Omega\sin(t\Omega) & \cos(t\Omega) \\
       \end{array}
     \right).
\end{equation}
Taking  this formula and the scheme of trigonometric integrators
\eqref{methods} into account, it is arrived that
\begin{eqnarray}
\left(
  \begin{array}{c}
d_{1} \\
\dot{d}_{1} \\
  \end{array}
\right) &=&  \kappa h\left(
          \begin{array}{c}
            h \bar{b}_{1}(V)f(u_{\frac{1}{2} }) \\
            b_{1}(V)f(u_{\frac{1}{2} }) \\
          \end{array}
        \right)
-\kappa\int_{0}^{h}R(h-t)\left(
                                 \begin{array}{c}
                                   0 \\
                                   f(u(\cdot,t)) \\
                                 \end{array}
                               \right)d t  \nonumber \\
   &= & \kappa h\left(
          \begin{array}{c}
            h \bar{b}_{1}(V)f(u_{\frac{1}{2} }) \\
            b_{1}(V)f(u_{\frac{1}{2} }) \\
          \end{array}
        \right)
- \kappa h R(\frac{h}{2})\left(
                   \begin{array}{c}
                     0 \\
                     f(u_{\frac{1}{2} }) \\
                   \end{array}
                 \right)
\label{sta er 1}\\
    &+ & \kappa h R(\frac{h}{2})\left(
                   \begin{array}{c}
                     0 \\
                     f(u_{\frac{1}{2} }) \\
                   \end{array}
                 \right)-\kappa h R(\frac{h}{2})\left(
                   \begin{array}{c}
                     0 \\
                       f(u(\cdot,\frac{h}{2})) \\
                   \end{array}
                 \right) \label{sta er 2}\\
    &+ &\kappa h R(\frac{h}{2})\left(
                   \begin{array}{c}
                     0 \\
                       f(u(\cdot,\frac{h}{2})) \\
                   \end{array}
                 \right)-\kappa\int_{0}^{h}R(h-t)\left(
                                 \begin{array}{c}
                                   0 \\
                                   f(u(\cdot,t)) \\
                                 \end{array}
                               \right)d t. \label{sta er 3}
\end{eqnarray}

In the following proof, we will use the   results (see
\cite{Gauckler17})
 \begin{equation}\label{dt pro}%
\normmm{\frac{d^l}{dt^l}R(t)\left(
                              \begin{array}{c}
                                u \\
                                \dot{v} \\
                              \end{array}
                            \right)
}_1=\normmm{(u,\dot{v} )}_{1+l}, \quad
\norm{\frac{d^{2-l}}{dt^{2-l}} f(u(\cdot,t)) }_{1+l}\leq C_M,
\end{equation}
where $l=0,1,2$.

   $\bullet$  Bound of \eqref{sta er 1}. According to  \eqref{matr R},    \eqref{sta er 1} is seen to be
   of the form
\begin{equation*}\begin{array}[c]{ll}
\kappa h\left(
          \begin{array}{c}
            h \bar{b}_{1}(V)f(u_{\frac{1}{2} }) \\
            \bar{b}_{1}(V)f(u_{\frac{1}{2} }) \\
          \end{array}
        \right)
- \kappa h R(\frac{h}{2})\left(
                   \begin{array}{c}
                     0 \\
                     f(u_{\frac{1}{2} }) \\
                   \end{array}
                 \right)
                 =\kappa h\left(
          \begin{array}{c}
            h \big(\bar{b}_{1}(V)-\frac{1}{2}\sinc(\frac{1}{2}h\Omega)\big)f(u_{\frac{1}{2} }) \\
            \big(b_{1}(V)- \cos(\frac{1}{2}h\Omega)\big)f(u_{\frac{1}{2} }) \\
          \end{array}
        \right).
\end{array}\end{equation*}
By the third and fifth formulae of \eqref{ass1}, we obtain
\begin{equation*}\begin{array}[c]{ll}
&\norm{\big(\bar{b}_{1}(V)-\frac{1}{2}\sinc(\frac{1}{2}h\Omega)\big)f(u_{\frac{1}{2}
})}_2\leq C \norm{ h\Omega f(u_{\frac{1}{2} })}_2\\
 =&C h \norm{ f(u_{\frac{1}{2} })}_3\leq C h
\Lambda_3(\norm{u}_{3})\norm{u}_{3}\norm{u}_{5}\leq C h,\\
&\norm{\big(b_{1}(V)- \cos(\frac{1}{2}h\Omega)\big)f(u_{\frac{1}{2}
})}_1\leq C \norm{ h^2\Omega^2 f(u_{\frac{1}{2} })}_1\\
=&C h^2 \norm{ f(u_{\frac{1}{2} })}_3\leq C h^2
\Lambda_3(\norm{u}_{3})\norm{u}_{3}\norm{u}_{5}\leq C h^2.
\end{array}\end{equation*}
Thus the term on right-hand side of \eqref{sta er 1} is bounded by
{$C h^3$.}

 $\bullet$  Bound of \eqref{sta er 2}.
For \eqref{sta er 2}, one has
\begin{equation*}\begin{array}[c]{ll}
&\normmm{\textmd{term  of }\eqref{sta er 2}}^2_1= \kappa^2
h^2\norm{\Omega^{-1}\sin(\frac{1}{2}h\Omega)(f(u_{\frac{1}{2}
})-f(u(\cdot,\frac{h}{2})))}_2^2\\
&\quad+\kappa^2 h^2 \norm{ \cos(\frac{1}{2}h\Omega)(f(u_{\frac{1}{2}
})-f(u(\cdot,\frac{h}{2})))}_1^2\\
&\leq \kappa^2 h^2\norm{f(u_{\frac{1}{2}
})-f(u(\cdot,\frac{h}{2}))}_1^2+\kappa^2 h^2\norm{ f(u_{\frac{1}{2}
})-f(u(\cdot,\frac{h}{2}))}_1^2.
\end{array}\end{equation*}
Since
\begin{equation*}\begin{array}[c]{ll}
&\norm{u_{\frac{1}{2} }-u(\cdot,\frac{h}{2})}_3=
\abs{\kappa}\norm{\int_{0}^{\frac{h}{2}}(\frac{h}{2}-t)\sinc((\frac{h}{2}-t)\Omega)f(u(\cdot,t))d
t}_3\\
\leq&\abs{\kappa}
\int_{0}^{\frac{h}{2}}(\frac{h}{2}-t)\norm{f(u(\cdot,t))}_3 d
t\leq\abs{\kappa}  C_M\int_{0}^{\frac{h}{2}}(\frac{h}{2}-t)  d t
\leq Ch^2\abs{\kappa},
\end{array}\end{equation*}
we obtain \begin{equation*}\begin{array}[c]{ll}
\norm{f(u_{\frac{1}{2} })-f(u(\cdot,\frac{h}{2}))}_1  &\leq
\Lambda_1\Big(\norm{u_{\frac{1}{2}}}_{3}+\norm{u(\cdot,\frac{h}{2}))}_{3}\Big)\Big(\norm{u_{\frac{1}{2}}}_{3}
+\norm{u(\cdot,\frac{h}{2}))}_{3}\Big)\\
&\norm{u_{\frac{1}{2}}-u(\cdot,\frac{h}{2}))}_{3}\leq
Ch^2\abs{\kappa},
\end{array}\end{equation*}
where \eqref{fu pro2} is used here.  Therefore, {it is arrived that
$\normmm{\textmd{term of }\eqref{sta er 2}}_1\leq C h^3. $}

 $\bullet$  Bound of \eqref{sta er 3}.
 The essential technology  used  here is the quadrature error of the mit-point rule. From its second-order
Peano kernel $K_2$, it follows that $ \textmd{term of } \eqref{sta
er 3}=-h^3 \kappa \int_{0}^{1}K_2(\sigma)l''(\sigma h)d\sigma $ with
$l(t)=R(h-t)\left(
                                 \begin{array}{c}
                                   0 \\
                                   f(u(\cdot,t)) \\
                                 \end{array}
                               \right).$ By \eqref{dt pro}, one arrives
$$\normmm{\textmd{term of } \eqref{sta er 3}}_1\leq C_M h^3
\abs{\kappa}. $$

All these estimates    together imply the result of this lemma. \QED

\begin{rem}
We remark that this lemma of the local error bound can be extended
to a nonzero $g$ in  \eqref{wave equa}  since the proof is only
based on the estimates \eqref{fu pro1} and \eqref{fu pro2}.
\end{rem}

\subsection{Proof of Theorem \ref{thm-error bounds}}

 \proof Denote by $C_1$ and $C_2$    the constants appearing in Propositions
\ref{pro-Stability} and \ref{prop-local error}, respectively.
 Let
{$h_0 = \sqrt{M/(C_2  T e^{C_1  T} )}$} and  it will be shown by
induction on $n$ that for {$h\leq h_0$}
\begin{equation}\label{er ma re 1}
\normmm{(u_n, \dot{u}_n)-(u(\cdot, t_n),\pa_t u(\cdot,
t_n))}_{1}\leq {C_2 e^{C_1nh}nh^3}\end{equation} as long as
$t_{n}=nh\leq T.$

 Firstly, it is obvious  that \eqref{er ma re 1} holds for $n=0.$  In what follows, it is assumed   that  \eqref{er ma re 1} holds for
$n=0,1,\ldots,m-1$. Choose $n=m-1$ and then we have
$$\normmm{(u_{m-1}, \dot{u}_{m-1})-(u(\cdot, t_{m-1}),\pa_t u(\cdot,
t_{m-1}))}_{1}\leq C_2\abs{\kappa}
e^{C_1\abs{\kappa}(m-1)h}(m-1)h^3,$$ which implies
$$\normmm{(u_{m-1}, \dot{u}_{m-1})}_{1}\leq M+{C_2
e^{C_1 (m-1)h}(m-1)h^3\leq2M}$$ as long as $t_{m-1}-t_0=(m-1)h\leq
T.$

For the global error, one has
\begin{equation*}
\begin{array}[c]{ll}
&\normmm{(u_{m}, \dot{u}_{m})-(u(\cdot, t_{m}),\pa_t u(\cdot,
t_{m}))}_{1}\\
 =& \normmm{\varphi_h(u_{m-1}, \dot{u}_{m-1})-(u(\cdot, t_{m}),\pa_t
u(\cdot,
t_{m}))}_{1} \nonumber \\
  \leq & \normmm{\varphi_h(u_{m-1}, \dot{u}_{m-1})-\varphi_h(u(\cdot,
t_{m-1}),\pa_t u(\cdot, t_{m-1}))}_{1}  \\
    &+  \normmm{\varphi_h(u(\cdot,
t_{m-1}),\pa_t u(\cdot, t_{m-1}))-(u(\cdot, t_{m}),\pa_t u(\cdot,
t_{m}))}_{1}.    \end{array}
\end{equation*}
From Proposition \ref{pro-Stability}, it follows that{
\begin{equation*}
\begin{array}[c]{ll}& \normmm{\varphi_h(u_{m-1}, \dot{u}_{m-1})-\varphi_h(u(\cdot,
t_{m-1}),\pa_t u(\cdot, t_{m-1}))}_{1} \\
\leq& (1+C_1  h) \normmm{(u_{m-1}, \dot{u}_{m-1})-(u(\cdot,
t_{m-1}),\pa_t u(\cdot, t_{m-1}))}_{1}
\\
\leq& (1+C_1  h)C_2 e^{C_1 (m-1)h}(m-1)h^3.
\end{array}
\end{equation*}}
On the other hand,  in the light of  Proposition  \ref{prop-local
error}, one reaches {
\begin{equation*}
\begin{array}[c]{ll}&\normmm{\varphi_h(u(\cdot,
t_{m-1}),\pa_t u(\cdot, t_{m-1}))-(u(\cdot, t_{m}),\pa_t u(\cdot,
t_{m}))}_{1} \leq C_2  h^{3}.
\end{array}
\end{equation*}}
 Thus,  it is obtained that {
\begin{equation*}
\begin{array}[c]{ll} & \normmm{(u_{m}, \dot{u}_{m})-(u(\cdot, t_{m}),\pa_t u(\cdot,
t_{m}))}_{1}\leq (1+C_1 h)C_2  e^{C_1 (m-1)h}(m-1)h^3+ C_2  h^{3}.
\end{array}
\end{equation*}}
Expanding  {$e^{C_1 (m-1)h}$} by   Taylor series, the right-hand
side of the above   inequality  becomes {\begin{equation*}
\begin{array}[c]{ll}  \ \ \ & (1+C_1 h)  C_2
e^{C_1(m-1)h}(m-1)h^3+C_2h^{3}\\
=& (1+C_1h)  C_2
\sum\limits_{k=0}^{\infty}\dfrac{(C_1(m-1)h)^{k}}{k!}
 (m-1)h^3+C_2h^{3}\\
 =&C_2\sum\limits_{k=0}^{\infty}\dfrac{1}{k!}(C_1)^k (m-1)^{k+1}h^{k+3}+
 C_2\sum\limits_{k=0}^{\infty}\dfrac{1}{k!}(C_1)^{k+1} (m-1)^{k+1}h^{k+4}+C_2h^{3}\\
  =&C_2mh^3 +C_2\sum\limits_{k=1}^{\infty}\big((m-1)^{k+1}+k(m-1)^{k}\big)\dfrac{1}{k!}(C_1)^k h^{k+3}.\\
\end{array}
\end{equation*}}
According to  the   fact
\begin{equation*}
\begin{array}[c]{ll}
(m-1)^{k+1}+k(m-1)^{k} \leq (m-1+1)^{k+1}=m^{k+1}\quad \textmd{for}
\quad m\geq1,\ \ k\geq1,
\end{array}
\end{equation*}
we obtain{
\begin{equation*}
\begin{array}[c]{ll}  \ \ \ & (1+C_1 h)  C_2
e^{C_1(m-1)h}(m-1)h^3+C_2h^{3}\\
\leq&C_2mh^3
+C_2\sum\limits_{k=1}^{\infty}m^{k+1}\dfrac{1}{k!}(C_1)^k h^{k+3}
=C_2 e^{C_1mh}mh^3.
\end{array}
\end{equation*}}
 Therefore, \eqref{er ma re 1} holds for $n=m.$
By  induction, it is true that{
\begin{equation*}
\normmm{(u_n, \dot{u}_n)-(u(\cdot, t_n),\pa_t u(\cdot,
t_n))}_{1}\leq C_2 e^{C_1T}T h^2\leq Ch^2,\end{equation*}} which
proves
 the statement of Theorem \ref{thm-error bounds}. \QED

\begin{rem}
The proof  also  holds for   a nonzero $g$ in  \eqref{wave equa}
since it is   based on Propositions \ref{pro-Stability} and
  \ref{prop-local error}, which are true for nonzero $g$.
\end{rem}

\section{Proof of error bounds for  full-discrete trigonometric integrators}\label{sec-Global error2}
In this section, we prove the error bounds for full-discrete
trigonometric integrators. Throughout the proof, we use  the
following approximation property of the $L^2$-orthogonal projection
$\mathcal{P}^{K}$:
\begin{equation}\label{K-2}
\norm{\mathcal{P}^{K}(v)}_{s}\leq \norm{v}_{s}\quad \textmd{for}
\quad v\in H^s
 \end{equation}
and
\begin{equation}\label{K-1}
\norm{v-\mathcal{P}^{K}(v)}_{s'}\leq K^{-(s-s')}\norm{v}_{s}\quad
\textmd{for} \quad v\in H^s,
 \end{equation}
where $s\geq s'\geq 0$. In addition, for the trigonometric
interpolation $\mathcal{I}^{K}$,
 we use the
approximation property for $s\geq s'\geq 0$ with  $s-
s'>\frac{1}{2}$
\begin{equation}\label{K-3}
\norm{v-\mathcal{I}^{K}(v)}_{s'}\leq C
_{s,s'}K^{-(s-s')}\norm{v}_{s}\quad \textmd{for} \quad v\in H^s
 \end{equation}
 and its stability
\begin{equation}\label{K-4}
\norm{\mathcal{I}^{K}(v)}_{s}\leq C _s\norm{v}_{s}\quad \textmd{for}
\quad v\in H^s.
 \end{equation}
It is noted that all estimates in the following are independent of
the spatial discretization parameter $K$.

\subsection{Stability}
The result of Lemma \ref{lem-remainder} can be extended to the
full-discrete trigonometric integrator  \eqref{methods-full}
directly.
\begin{lem}\label{lem-remainder-full}
Under the conditions given in Lemma \ref{lem-remainder}, it follows
that
\begin{equation*}\begin{array}[c]{ll}
&\normmm{(u^K_{n+1}-v^K_{n+1},\dot{u}^K_{n+1}-\dot{v}^K_{n+1})}^2_{1}\\
=&\normmm{(u^K_{n}-v^K_{n},\dot{u}^K_{n}-\dot{v}^K_{n})}^2_{1}+\kappa
\mathcal{R}^K(u^K_{n+1},v^K_{n+1},u^K_{n},v^K_{n})\end{array}\end{equation*}
with the remainder
\begin{equation}\begin{array}[c]{ll}\label{remainder-full}
&\mathcal{R}^K(u^K_{n+1},v^K_{n+1},u^K_{n},v^K_{n}) \\
=& \langle
2\sinc(h\Omega)^{-1}\bar{b}_{1}(V)(u^K_{n+1}-u^K_{n}-v^K_{n+1}+v^K_{n}),\hat{f}^K(u^K_{n+\frac{1}{2}})-\hat{f}^K(v^K_{n+\frac{1}{2}})\rangle_1.\end{array}\end{equation}
\end{lem}

For this remainder $\mathcal{R}^K$, we have the following bound.
\begin{lem}\label{lem-remainder bound-full}
(\textbf{Bound of the remainder.}) Under the conditions given in
Lemma \ref{lem-remainder bound},  the remainder $\mathcal{R}^K$ is
bounded by
\begin{equation}\begin{array}[c]{ll}\label{remainder bound-full}
\abs{{\kappa}\mathcal{R}^K(u^K_{n+1},v^K_{n+1},u^K_{n},v^K_{n})
}\leq C_M
h\normmm{(u^K_{n}-v^K_{n},\dot{u}^K_{n}-\dot{v}^K_{n})}_{1}^2.\end{array}\end{equation}
\end{lem}
\proof This lemma is proved in a similar way to that   of Lemma
\ref{lem-remainder bound} by using in addition  the bounds
\eqref{K-2} and \eqref{K-4} on $\mathcal{P}^{K}$ and
$\mathcal{I}^{K}$ and  the property
\begin{equation}\label{add-pro}
\langle v^K,\mathcal{P}^{K}(w)\rangle_s=\langle v^K, w\rangle_s
\quad for \quad v^K\in \mathcal{V}^K,\ w\in H^s\end{equation}
  with $s = 1$.  \QED

The stability of the full-discrete trigonometric integrator
\eqref{methods-full} is obtained immediately by these two lemmas.
\begin{prop}\label{pro-Stability-full}
(\textbf{Stability.}) Under the conditions given in Proposition
\ref{pro-Stability-full},   we have
\begin{equation}\begin{array}[c]{ll}\label{Stability-full}
\normmm{(u^K_{n+1},
\dot{u}^K_{n+1})-\varphi^K_h(u^K(\cdot,t_n),\pa_t
u^K(\cdot,t_n))}_1^2\leq {(1+C_M h)}
\normmm{(e^K_n,\dot{e}^K_n)}_1^2,
\end{array}\end{equation}
where the global error $(e^K_n,\dot{e}^K_n)$ is defined by
$$(e^K_n,\dot{e}^K_n)=(u^K_{n}, \dot{u}^K_{n})-(u^K(\cdot,t_n),\pa_t
u^K(\cdot,t_n)).$$ \end{prop}

\subsection{Local error bound}
For full-discrete trigonometric integrator  \eqref{methods-full},
the local error bound is presented as follows.

\begin{prop}\label{prop-local error-full}
(\textbf{Local error in $H^{2}\times H^{1}$.}) Under the conditions
of Proposition \ref{prop-local error}, one has
$$\normmm{(d^K_{n+1},\dot{d}^K_{n+1})}_{1}\leq {C_Mh^{3} +C_Mh K^{-s-2}},$$
where the local error $(d^K_{n+1},\dot{d}^K_{n+1})$ is defined by
$$(d^K_{n+1},\dot{d}^K_{n+1})=\varphi^K_h(u^K(\cdot,t_n),\pa_t
u^K(\cdot,t_n))-(u^K(\cdot,t_{n+1}),\pa_t u^K(\cdot,t_{n+1})).$$
\end{prop}
\proof Similar to the proof of Proposition \ref{prop-local error},
we only consider the case $n = 0$. Using this formula and letting
$\tilde{f}^{K}(u)= \mathcal{P}^{K} \circ f$, the local error can be
rewritten in the form
\begin{eqnarray}
\left(
  \begin{array}{c}
d^K_{1} \\
\dot{d}^K_{1} \\
  \end{array}
\right) &=&  \kappa h\left(
          \begin{array}{c}
            h \bar{b}_{1}(V)\hat{f}^K(u^K_{\frac{1}{2} }) \\
            b_{1}(V)\hat{f}^K(u^K_{\frac{1}{2} }) \\
          \end{array}
        \right)
-\kappa\int_{0}^{h}R(h-t)\left(
                                 \begin{array}{c}
                                   0 \\
                                   \tilde{f}^{K}(u(\cdot,t)) \\
                                 \end{array}
                               \right)d t  \nonumber \\
   &= & \kappa h\left(
          \begin{array}{c}
            h \bar{b}_{1}(V)\hat{f}^K(u^K_{\frac{1}{2} }) \\
            b_{1}(V)\hat{f}^K(u^K_{\frac{1}{2} }) \\
          \end{array}
        \right)
- \kappa h R(\frac{h}{2})\left(
                   \begin{array}{c}
                     0 \\
                     \hat{f}^K(u^K_{\frac{1}{2} }) \\
                   \end{array}
                 \right)
\label{sta er 1-full}\\
    &+ & \kappa h R(\frac{h}{2})\left(
                   \begin{array}{c}
                     0 \\
                     \hat{f}^K(u^K_{\frac{1}{2} }) \\
                   \end{array}
                 \right)-\kappa h R(\frac{h}{2})\left(
                   \begin{array}{c}
                     0 \\
                       \tilde{f}^{K}(u^K_{\frac{1}{2} }) \\
                   \end{array}
                 \right) \label{sta er 2-full}\\
 &+ &\kappa h R(\frac{h}{2})\left(
                   \begin{array}{c}
                     0 \\
                       \tilde{f}^{K}(u^K_{\frac{1}{2} })  \\
                   \end{array}
                 \right)-\kappa h R(\frac{h}{2})\left(
                   \begin{array}{c}
                     0 \\
                       \tilde{f}^{K}(u(\cdot,\frac{h}{2})) \\
                   \end{array}
                 \right)\label{sta er 3-full}\\
    &+ &\kappa h R(\frac{h}{2})\left(
                   \begin{array}{c}
                     0 \\
                       \tilde{f}^{K}(u(\cdot,\frac{h}{2})) \\
                   \end{array}
                 \right)-\kappa\int_{0}^{h}R(h-t)\left(
                                 \begin{array}{c}
                                   0 \\
                                   \tilde{f}^{K}(u(\cdot,t)) \\
                                 \end{array}
                               \right)d t. \label{sta er 4-full}
\end{eqnarray}
 Bounds of \eqref{sta er 1-full}, \eqref{sta er 3-full} and \eqref{sta er 4-full}  can be derived by using the same way as
in the proof of Proposition \ref{prop-local error} and by using in
addition the properties \eqref{K-1}-\eqref{K-4} of $\mathcal{P}^K $
and $\mathcal{I}^K$  and the assumed regularity of the exact
solution.  For the bound of \eqref{sta er 2-full}, we have
$$\hat{f}^K-\tilde{f}^{K}=\mathcal{P}^{K}
\circ (f^K-f).$$ By the   arguments of the proof of Proposition
\ref{prop-local error} as well as \eqref{K-2} and \eqref{K-3}, the
estimate $C_M h K^{-s-4}\abs{\kappa}$ in $H^2 \times H^1$ and $C_M h
K^{-s-3}\abs{\kappa}$ in $H^3 \times H^2$ for \eqref{sta er 2-full}
can be obtained.
  \QED

\subsection{Proof of Theorem \ref{thm-error bounds full}}
\proof Based on the above analysis given in this section, the proof
of Theorem \ref{thm-error bounds full} is similar to that of Theorem
\ref{thm-error bounds} with some obvious adjustments. \QED

\begin{rem}
It is noted that the proof of error bounds  for full-discrete
trigonometric integrators  does not require any CFL-type coupling of
the discretization parameters.
\end{rem}
\section{Proof of Theorem \ref{thm-error bounds a}} \label{sec-simple proof}
We  consider the following Strang splitting method
 \begin{equation*}
\begin{array}[c]{ll}
&\textmd{1. } (q_{+}^{n},p_{+}^{n}) = \Phi_{h/2,\textmd{L}}
(q^{n},p^{n}): \\
&\left(
             \begin{array}{c}
              q_{+}^{n} \\
              p_{+}^{n} \\
             \end{array}
           \right)=\left(
                                                           \begin{array}{cc}
                                                            \cos(\frac{h\Omega}{2})& \Omega^{-1}\sin(\frac{h\Omega}{2}) \\
                                                             -\Omega \sin(\frac{h\Omega}{2}) & \cos(\frac{h\Omega}{2}) \\
                                                           \end{array}
                                                         \right)\left(
             \begin{array}{c}
               q^{n} \\
               p^{n} \\
             \end{array}
           \right),\\
&\textmd{2. }(q_{-}^{n},p_{-}^{n}) = \Phi_{h,\textmd{NL}}
(q_{+}^{n},p_{+}^{n}):\\
&\left(
             \begin{array}{c}
              q_{-}^{n} \\
              p_{-}^{n} \\
             \end{array}
           \right)=\left(
             \begin{array}{c}
              q_{+}^{n} \\
              p_{+}^{n}+h\Upsilon(h \Omega)g(q_{-}^{n}) \\
             \end{array}
           \right),\\
&\textmd{3. }(q^{n+1},p^{n+1}) = \Phi_{h/2,\textmd{L}}
(q_{-}^{n},p_{-}^{n}): \\
&\left(
             \begin{array}{c}
              q^{n+1} \\
              p^{n+1} \\
             \end{array}
           \right)=\left(
                                                           \begin{array}{cc}
                                                            \cos(\frac{h\Omega}{2})& \Omega^{-1}\sin(\frac{h\Omega}{2}) \\
                                                             -\Omega \sin(\frac{h\Omega}{2}) & \cos(\frac{h\Omega}{2}) \\
                                                           \end{array}
                                                         \right)\left(
             \begin{array}{c}
               q_{-}^{n} \\
               p_{-}^{n}\\
             \end{array}
           \right).
\end{array}
\end{equation*}
It can be checked that the trigonometric integrator $\phi_h$ can be
expressed by this Strang splitting as
 \begin{equation}
\phi_h=\Phi_{h/2,\textmd{L}} \circ \Phi_{h,\textmd{NL}} \circ
\Phi_{h/2,\textmd{L}}
  \label{connetciton1}%
\end{equation}
if and only if
\begin{equation}  \label{Upsilon}
\Upsilon(h \Omega)=b_1(h \Omega)\cos^{-1}\big(\frac{1}{2}h
\Omega\big)=2\bar{b}_1(h \Omega)\textmd{sinc}^{-1}\big(\frac{1}{2}h
\Omega\big).\end{equation}

 On the other hand, for the Strang splitting
 \begin{equation*}
\hat{\phi}_h=\Phi_{h/2,\textmd{NL}} \circ \Phi_{h,\textmd{L}} \circ
\Phi_{h/2,\textmd{NL}},
\end{equation*}
it is identical to  a    trigonometric integrator
((XIII.2.7)--(XIII.2.8) given on p.481 of \cite{hairer2002})
\begin{equation}\label{TI}
\begin{array}[c]{ll} &q^{n+1}=
\cos(h\Omega)q^{n}+h\textmd{sinc}(h\Omega)p^{n}+\frac{1}{2}h^2
 \textmd{sinc}(h\Omega)\Upsilon(h
\Omega)g(q^n),\\
 &p^{n+1}=-\Omega
 \sin(h\Omega)q^{n}+\cos(h\Omega)p^{n}+\frac{1}{2}h\big(
 \cos(h\Omega)\Upsilon(h \Omega)g(q^n)+\Upsilon(h
\Omega)g(q^{n+1})\big).
\end{array}
\end{equation}
This  trigonometric integrator with the choice
\begin{equation}\label{scc}\Upsilon(h
\Omega)=\textmd{sinc}(h\Omega)\end{equation} has been discussed in
\cite{Gauckler17} for quasilinear wave equations.
 Thus based on the following important connection
 \begin{equation}\begin{aligned}
\underbrace{\phi_h\circ \cdots \circ\phi_h}_{n\ \textmd{times}}
=&\Phi_{h/2,\textmd{L}}\circ\Phi_{h/2,\textmd{NL}}\circ\big(\underbrace{\hat{\phi}_h\circ
\cdots \circ \hat{\phi}_h}_{n-1\ \textmd{times}}\big)   \circ
\Phi_{h/2,\textmd{NL}}\circ \Phi_{h/2,\textmd{L}}.
  \label{connetciton3}%
\end{aligned}\end{equation}
it is arrived that TI2 and the trigonometric integrator
\eqref{TI}-\eqref{scc} have similar error bounds when they are used
to solving quasilinear wave equations. Therefore, the error bounds
of TI2 are obtained immediately by considering the results given in
\cite{Gauckler17} and the proof of Theorem \ref{thm-error bounds a}
is complete.
\section{Concluding remarks}\label{sec-conclu}
This paper studied  error bounds of  one-stage explicit
 trigonometric integrators for solving quasilinear wave equations. Second-order convergence for the
semidiscretization in time was proved and the error bounds of  fully
discrete scheme were also presented without requiring any CFL-type
coupling of the discretization parameters.

Last but not least, the analysis of trigonometric integrators in
this paper can be extended to quasilinear wave equations \eqref{wave
equa} without Klein-Gordon term $-u$  and also works for higher
spatial dimensions. The application and analysis of  trigonometric
integrators for quasilinear wave equations with $\kappa=1$ and for
more general quasilinear wave equations or other kinds of PDEs will
be our future work.

\section*{Acknowledgements}

The first author is grateful to Ludwig Gauckler for pointing out the
connection between symmetric trigonometric methods given in this
paper and symmetric trigonometric integrators of \cite{hairer2002},
which motives Section \ref{sec-simple proof}. We also thank him for
the very helpful comments on the proof of Lemma \ref{lem-remainder
bound}.

\end{document}